\documentclass[11pt]{article}
\usepackage[letterpaper,margin=1in]{geometry}
\usepackage[T1]{fontenc}
\usepackage{lmodern}
\usepackage{microtype}
\usepackage{amsmath,amssymb,amsthm,mathtools}
\usepackage{booktabs,caption,needspace,setspace}
\usepackage[numbers,square,sort&compress]{natbib}
\usepackage[colorlinks=true,linkcolor=black,citecolor=black,urlcolor=black]{hyperref}
\ifdefined\pdfgentounicode
  \input{glyphtounicode}
\fi
\hypersetup{bookmarksnumbered=true,bookmarksdepth=2,
  pdftitle={Unimodality of generalized subword counts in one-descent and bitonic compositions},
  pdfauthor={Yibo Zhang and Zhanhe Zhang},
  pdfkeywords={integer compositions, generalized subword order, rank-unimodality, layered permutations}}
\numberwithin{equation}{section}
\newtheorem{theorem}{Theorem}[section]
\newtheorem{lemma}[theorem]{Lemma}
\newtheorem{proposition}[theorem]{Proposition}
\newtheorem{corollary}[theorem]{Corollary}
\theoremstyle{definition}
\newtheorem{condition}[theorem]{Condition}

\newcommand{\wt}{\operatorname{wt}}
\newcommand{\coef}[2]{[z^{#1}]#2}
\newcommand{\Prb}{\mathbb P}
\newcommand{\firstmode}{m_-}
\newcommand{\lastmode}{m_+}
\newcommand{\orcidlink}[1]{\href{https://orcid.org/#1}{\nolinkurl{#1}}}

\title{Unimodality of generalized subword counts\\
in one-descent and bitonic compositions}
\author{Yibo Zhang$^a$ {\small(\orcidlink{0009-0002-2005-5781})},
  Zhanhe Zhang$^a$ {\small(\orcidlink{0009-0007-2883-9640})}\\[0.6em]
  \small $^a$Central University of Finance and Economics, Beijing, China\\[0.4em]
  \small Both authors contributed equally and share first authorship.\\
  \small Correspondence to Zhanhe Zhang at
  \href{mailto:2024311337@email.cufe.edu.cn}{2024311337@email.cufe.edu.cn}.}
\date{November 25, 2026}

\begin{document}
\maketitle
\begin{abstract}
Vatter conjectured that the numbers of distinct generalized subwords of
successive weights in a composition form a unimodal sequence. We prove this
for compositions with at most one strict decrease between adjacent parts
and for bitonic compositions, whose parts first weakly increase and then
weakly decrease. Both families reduce to a peak larger than every part in
its two monotone arms. The rank polynomial at the peak involves subtracting
the product of the two arm rank polynomials.
A coupling of uniform subwords at consecutive weight bounds relates
the positive rank differences in a single arm to its cumulative counts
at earlier ranks. Let \(p_k\) be the coefficients of the reversed product
and \(u_k=\max(p_k-p_{k-1},0)\), with \(p_{-1}=0\).
We prove that \(2^k u_k\) is unimodal. From its first maximum onward,
each \(u_k\) bounds the sum of all later positive differences and thus
controls the subtraction at the peak. The proof uses exact finite
enumeration for one local inequality when some composition of weight
at most twelve is missing from one arm.
\end{abstract}
\noindent\textit{Keywords:} Integer compositions; generalized subword order; rank-unimodality; layered permutations.

\noindent\textit{MSC 2020:} 05A20 (primary), 06A07, 05A15.

\onehalfspacing
\clearpage

\section{Introduction}\label{sec:introduction}

A composition \(w=(w_1,\ldots,w_m)\) is a finite sequence of positive
integers, with the empty composition allowed. Its weight \(\wt(w)\)
is the sum of its parts. In generalized subword order \citep{mcnamara-sagan},
we write \(u=(u_1,\ldots,u_k)\preceq w\) if it has an embedding in \(w\),
that is, indices \(1\le i_1<\cdots<i_k\le m\) with \(u_j\le w_{i_j}\)
for every \(j\).
The principal ideal of \(w\) consists of the distinct \(u\preceq w\),
ranked by weight. Let \(a_n(w)\) be the number at rank \(n\), including
\(a_0(w)=1\) for the empty subword.

A sequence is unimodal if it is weakly increasing up to a maximum and
weakly decreasing thereafter. The ideal is rank-unimodal if
\(a_0(w),\ldots,a_{\wt(w)}(w)\) is unimodal.
\citet{vatter} asked whether this holds for every composition
(Conjecture~2.6 in \href{https://arxiv.org/abs/2602.16355v3}{arXiv:2602.16355v3}).
A descent is a strict decrease between adjacent parts.
A composition is bitonic if its parts weakly increase to some position
and then weakly decrease. We include the empty composition in this class.

\begin{theorem}\label{thm:main}
If a composition \(w\) has at most one descent or is bitonic, then
\(a_0(w),\ldots,a_{\wt(w)}(w)\) is unimodal.
\end{theorem}

A layered permutation has decreasing blocks of sizes \(w_1,\ldots,w_m\),
with smaller values in each block than in the next. The positive sizes
retained from the original blocks are the lengths of the pattern's
maximal decreasing runs. Thus \(a_n(w)\) counts the distinct patterns
of length \(n\) in the corresponding layered permutation
\citep[Section~7.2]{sagan-vatter}, whose principal ideal is rank-unimodal
under either hypothesis of Theorem~\ref{thm:main}.
\citet[Conjecture~9.4]{mcnamara-steingrimsson} asked for unimodality
in every interval of the classical pattern poset. In consecutive
pattern order, where the selected positions must be a contiguous interval,
\citet[Corollary~5.2]{elizalde-mcnamara} proved that every interval is
rank-unimodal.

When \(w\) is weakly decreasing, every subword of length \(k\) fits
in the first \(k\) positions. Its rank sizes thus agree with those
in componentwise order, where \citet[Theorem~3.3(b)]{sagan} proved
unimodality. Reversal gives the weakly increasing case
\citep[Corollary~2.5]{vatter}. For one descent,
Zhanhe~Zhang (\href{https://doi.org/10.2139/ssrn.7397439}{SSRN~7397439})
split the rank polynomial at a maximum, proving unimodality
if a single part follows the descent or if the peak is at least the
combined weight of the arms minus two. The same paper covers one-descent
hosts with at most four parts.
Lizhong~Chen (\href{https://doi.org/10.2139/ssrn.7385218}{SSRN~7385218})
proved unimodality for every host with at most six parts.
For bitonic hosts Chen also proved log-concavity from the first strict
descent of the rank sequence when the parts are bounded by five or
one arm is a single part or a string of ones. Theorem~\ref{thm:main}
allows arbitrary lengths and part sizes in both families.

Write \(F_w(z)=\sum_{n=0}^{\wt(w)}a_n(w)z^n\) for the rank polynomial
and \(S_t=z+\cdots+z^t\) for a positive integer \(t\).
Consider the one-descent composition \((1,4,2,3,5)\). Removing the
final \(5\) leaves \((1,4,2,3)\), which has the same rank polynomial as
the bitonic word \((1,4,3,2)\). Splitting at \(4\) gives this equality
because reversal preserves the right arm's rank polynomial. In general,
let \(M\) be the part just before the descent. Remove the suffix of parts
at least \(M\) after the descent and keep one copy of \(M\) at the descent.
For a bitonic word, keep one maximum and reverse the decreasing arm.
Restoring a removed part of size \(t\) replaces the current rank polynomial
\(F\) by \(1+S_tF\), preserving unimodality. Section~\ref{sec:peaks}
justifies these reductions.

The remaining peak \(M\) exceeds every part in either arm. If the arm
rank polynomials are \(A,B\), the decomposition gives
\(A+B+(S_M-1)AB\). Let \(p_k\) be the coefficient at index \(k\) in the
coefficient reversal of \(AB\), and let \(u_k=\max(p_k-p_{k-1},0)\),
with \(p_{-1}=0\). The coupling in Theorem~\ref{thm:coupling} and the
single-arm comparison in Corollary~\ref{cor:arm} bound the growth of
positive differences in one arm. Proposition~\ref{prop:count-induction}
uses that comparison to prove unimodality of \(2^k u_k\). From its first
maximum onward, \(u_{k+1}\le u_k/2\), so the gain \(p_{k+M}-p_k\) is at most
\(u_k\). This is the bound needed to handle the subtraction at the peak.
\section{Preliminaries and rank polynomials}\label{sec:preliminaries}

The weight function grades generalized subword order. Given an embedding
of \(u\) in \(w\), increase the matched parts one unit at a time to their
values in \(w\). Insert a one at each unused position and increase it
to the required value. Every intermediate word lies below \(w\), and
each step increases the weight by one. Thus every cover raises the
weight by one. For nonempty \(w\), the coefficients at ranks
\(0,1,\wt(w)\) are one, corresponding to the empty word, \((1)\),
and \(w\), respectively. The maximum part of the empty word is zero
by convention.

For a polynomial \(V\) we write \(V_i=\coef{i}{V}\), with coefficients
outside its support taken to be zero. A nonnegative sequence \(x\) is
log-concave if \(x_i^2\ge x_{i-1}x_{i+1}\) at every interior index.
In ordinary subword order the selected letters retain their values.
\citet{chase} proved that the numbers of distinct subwords of successive
lengths in a fixed word form a log-concave sequence. For a nonzero
finitely supported unimodal sequence \(x\), let \(\firstmode(x)\) and
\(\lastmode(x)\) be the first and last indices of its maximum.

For a composition \(w\), use the cumulative counts
\[
 G_n(w)=\sum_{j\le n}a_j(w).
\]
In particular \(G_n(w)=0\) for \(n<0\) and is constant for \(n\ge\wt(w)\).
Write \(\Delta x_i=x_i-x_{i-1}\) and \(y_+=\max(y,0)\).
For forward rank sequences put \(d_n(w)=(\Delta a_n(w))_+\) for \(n\ge2\)
and \(d_n(w)=0\) for \(n\le1\). For reversed polynomials we take
positive differences from rank zero, using the zero coefficient at
rank \(-1\).

Let \([t]=1+z+\cdots+z^{t-1}\), \(S_t=z[t]\) for \(t\ge1\), and
\([0]=S_0=0\). Reversing the coefficients of \(F_w\) gives
\(q_w(z)=z^{\wt(w)}F_w(1/z)\). Reversing the \emph{parts} of \(w\)
preserves \(F_w\), since it reverses each embedded word. An increasing
arm is a composition with weakly increasing parts, possibly empty.
An \emph{interval arm} is an increasing arm \(X\) for which
\(F_X=[\wt(X)+1]\). Juxtaposition denotes concatenation of compositions
or parts.

\begin{lemma}[Rightmost placements and maximum decomposition]\label{lem:count}
Let \(X,Y\) be compositions. For \(t\ge\max(X)\) and \(t\ge1\),
the append formula is
\begin{equation}\label{eq:append}
 F_{Xt}=1+S_tF_X,\qquad q_{Xt}=[t]q_X+z^{\wt(X)+t}.
\end{equation}
For \(M\ge1\) and \(M\ge\max(X),\max(Y)\), put \(A=F_X\), \(B=F_Y\). Then
\begin{equation}\label{eq:core}
 F_{XMY}=A+B+(S_M-1)AB.
\end{equation}
Both statements allow empty arms.
\end{lemma}
\begin{proof}
Every nonempty subword of \(Xt\) has a unique expression \(uj\), with
\(u\preceq X\) and \(1\le j\le t\). Indeed, the last letter of any
embedding may be moved to the final position of \(Xt\). Counting by \(j\)
gives the first identity of \eqref{eq:append}, and reversing its
coefficients gives the second.

A subword of \(XMY\) that does not fit in \(X\) has a unique longest
prefix \(u\) that does. Its next letter \(j\) must lie at or beyond the
displayed \(M\) in every embedding. We can place \(j\) at \(M\), followed
by an arbitrary subword of \(Y\). Among all extensions \(uj\) with
\(u\preceq X\) and \(1\le j\le M\), those still in \(X\) are precisely
its nonempty subwords. The others have polynomial \(S_M A-(A-1)\).
Adding \(A\) to the product of this expression with \(B\) gives
\eqref{eq:core}.
\end{proof}

For a nonempty increasing arm we may replace its first part \(c>1\)
with \(c\) ones. Both initial words have polynomial \([c+1]\), and every
subsequent part is at least \(c\), so \eqref{eq:append} preserves the
rank polynomial at each step. This normalization preserves weight and
gives an arm \(1^\alpha(t_1,\ldots,t_h)\), where \(1^\alpha\) denotes
\(\alpha\ge1\) copies of one and \(2\le t_1\le\cdots\le t_h\) if \(h>0\).

\begin{lemma}[Elementary bounds]\label{lem:count-bounds}
For every host \(w\) of weight \(N\), \(a_n(w)>0\) for \(0\le n\le N\).
For \(n\ge1\),
\[
 a_n(w)\le 2a_{n-1}(w), \qquad
 a_n(w)\le\sum_{i<n}a_i(w).
\]
The sequence is nondecreasing through rank \(\lceil N/2\rceil\).
\end{lemma}
\begin{proof}
Lower the last part by one, deleting it if it equals one. Iteration from
\(w\) visits every weight from \(N\) to zero. A word has at most two
preimages under this operation, obtained by increasing its last part or
appending a one. This gives the factor of two. Deleting the entire last
part instead gives an injection from rank \(n\) into the union of the
earlier ranks, since \(n\) determines the missing part.

For the last assertion, equip the entries of the layered permutation of
\(w\) with their position order and their value order. The isomorphism
types induced on \(r\) entries are exactly the patterns of length \(r\).
This equivalence on subsets is hereditary, since an isomorphism pairs
the subsets of each smaller size pattern by pattern. The rank inequalities
of \citet[Corollary~2.11]{pouzet} give nondecrease through
\(\lfloor N/2\rfloor\) and
\(a_r(w)\le a_{N-r}(w)\) for \(2r\le N\). For odd \(N\), the latter
inequality gives the middle step.
\end{proof}

For a unimodal polynomial with nonnegative integer coefficients, the
indices at which \(a_i\ge h\) form an interval \([l_h,r_h]\) for each
positive integer \(h\) up to the largest coefficient. Summing the
interval polynomials \(\sum_{j=l_h}^{r_h}z^j\) recovers the polynomial.
The product of two interval polynomials counts their overlap as one
interval moves past the other. It is symmetric and unimodal about
the sum of their centers. Its height intervals have this common center,
so the same conclusion follows inductively for any number of factors.

For an integer \(c\), all height intervals have endpoint sum at least
\(c\) if and only if \(a_i\le a_j\) for all integers \(i\le j\) with
\(i+j\le c\). In one direction, each interval containing \(i\) must
then contain \(j\). Conversely, an interval \([l,r]\) with \(l+r<c\)
contains \(l\) and omits its reflection \(c-l\), so the coefficient
at \(l\) exceeds the one at \(c-l\).
Likewise, endpoint sums at most \(c\) give \(a_i\ge a_j\) for
\(i\le j\) with \(i+j\ge c\). Reversal \(i\mapsto d-i\) sends
an endpoint sum \(c\) to \(2d-c\).

\begin{lemma}[Uniform windows and arm layers]\label{lem:layers}
If \(P\) is nonnegative and unimodal, then \([t]P\) is unimodal for
every \(t\ge1\). For an increasing arm \(X\) of weight \(a\), \(F_X\)
is unimodal. Its lowest height interval is \([0,a]\), and all other
height intervals have endpoint sum at least \(a+1\).
If \(F_X\) is more than an interval, its second height interval is
\([2,a-1]\).
\end{lemma}
\begin{proof}
The coefficients of \([t]\) are log-concave, so discrete strong
unimodality gives the first assertion \citep{keilson-gerber}. One can
also see this from the differences \(P_n-P_{n-t}\), extending \(P\)
by zero. They are nonnegative through the last mode of \(P\).
Thereafter the entering coefficient decreases weakly and the departing
one increases weakly until its index reaches the first mode. From
then on \(P_n\le P_{n-t}\).

Append \(t\) to an increasing arm \(X\) of weight \(a\).
The window \(S_tF_X\) begins with coefficient one at rank one, so
adding its new constant term preserves unimodality. By induction, the
interval products in \(S_tF_X\) have endpoint sums at least \(a+t+1\).
Hence \((S_tF_X)_i\le(S_tF_X)_j\) for \(i\le j\) and \(i+j\le a+t+1\).
Adding the constant extends its lowest height interval to \([0,a+t]\)
and leaves the higher intervals unchanged. If \(Xt\) is noninterval,
\(X\) is nonempty and \(t\ge2\). The product of its lowest prefix
interval \([0,a]\) and \(S_t\) contributes at least two at every rank
from \(2\) through \(a+t-1\). Ranks zero, one, and \(a+t\) have
coefficient one, so the second height interval is exactly \([2,a+t-1]\).
The empty arm, a single part, and a string of ones each have an interval
as their rank polynomial.
\end{proof}
\section{A coupling at consecutive weight bounds}\label{sec:transport}

We call a host monotone if its parts are weakly increasing or weakly
decreasing. For such a host, increasing the weight bound from \(n\) to
\(n+1\) changes the distribution of both the weight and the length of
a uniformly chosen subword. We couple the two choices so that the weight
increases by at most one and the length never decreases. Reversal
preserves both statistics. We may therefore take a decreasing host, in
which every subword of length \(\ell\) fits in the first \(\ell\) positions.

Nonnegative sequences \(x,y\) are in ratio order if
\(x_i y_j\ge x_j y_i\) for all \(i<j\). If \(x\) is log-concave with
interval support, comparison of consecutive ratios gives this order for
every right translate
\(y_i=x_{i-s}\), \(s\ge0\). When \(x_jy_i>0\), the intermediate terms are positive by interval
support. When \(x_jy_i=0\), the cross-product inequality is immediate.
Linearity extends it to any nonnegative sum of right translates.

Equivalently, the Toeplitz matrix \((x_{j-i})\) of a nonnegative sequence
with interval support has nonnegative minors of order two precisely when
\(x\) is log-concave. This follows from the same ratio comparison, with
adjacent rows and columns giving the converse. The convolution of two
finite log-concave sequences with interval support is again log-concave
\citep{hoggar}. Indeed, the Toeplitz matrix of the convolution is the
product of the two Toeplitz matrices. Its minors of order two are
nonnegative by Cauchy--Binet, and the support is the sum of two
intervals. The same argument with the two-row matrix whose rows are
\(x,y\) shows that convolving both with a fixed log-concave sequence
of interval support preserves their ratio order. Only finitely many
indices contribute to each minor here, even for convolution with
\((1,1,\ldots)\), since the sequences are supported on a half-line.

For a decreasing host \(b_0\ge\cdots\ge b_{H-1}\ge1\), put \(b_H=0\)
and
\[
 P_0=1,\qquad P_\ell=\prod_{r=0}^{\ell-1}S_{b_r}
 \quad(1\le\ell\le H),\qquad p_\ell(k)=\coef{k}{P_\ell}.
\]
Thus \(P_\ell\) counts the subwords of length \(\ell\), and
\(F=\sum_{\ell=0}^H P_\ell\) counts the whole ideal. Write
\(C_\ell(n)=\sum_{k\le n}p_\ell(k)\). Their sum
\(C(n)=\sum_\ell C_\ell(n)\) is the cumulative count \(G_n\) of this host.

\begin{theorem}[Adjacent budgets]\label{thm:coupling}
For a monotone host and a budget \(n\ge0\), choose a word uniformly
among the subwords of weight at most \(n\). Let \(K_n\) be its
weight and \(J_n\) its length. The choices at budgets \(n,n+1\)
admit a coupling in which \(K_n,J_n,n-K_n\) all increase weakly.
For the length counts \(C_\ell(n)\) above, the cumulative count \(C(n)\)
is log-concave. For any nonincreasing sequence of nonnegative integer
thresholds \(s_0,\ldots,s_H\),
\begin{equation}\label{eq:threshold}
 \frac{\sum_{\ell=0}^H C_\ell(n-s_\ell)}{C(n)}
 \quad\text{is nondecreasing in }n\ge0.
\end{equation}
\end{theorem}
\begin{proof}
Condition first on the weight. If \(i<j\), the polynomial \(P_j\)
is \(P_i\) times a product of interval factors. The coefficients of
\(P_i\) are log-concave, so comparison with its translates gives, for \(k<k'\),
\begin{equation}\label{eq:tp}
 p_i(k)p_j(k')\ge p_i(k')p_j(k).
\end{equation}
Sum \eqref{eq:tp} over \(i<d\le j\) and divide by the two rank
counts at \(k,k'\), taking weights that occur. This gives
\(\Prb(J\ge d\mid K=k)\le\Prb(J\ge d\mid K=k')\).
Sampling the two conditional lengths at a common quantile makes the length at the higher weight at least as large.

Suppose the ratios \(C(i)/C(i-1)\) are nonincreasing for
\(1\le i\le n\), where \(n\ge1\). We construct the coupling from
budget \(n-1\) to budget \(n\). The weight distribution function
at budget \(v\) is \(C(\min(k,v))/C(v)\) at an integer \(k\).
Increasing the budget stochastically increases \(K\).
The same quantile increases it by at most one, because
\[
 \frac{C(k)}{C(n-1)}
 \le \frac{C(k+1)}{C(n)}
 \qquad(0\le k\le n-1).
\]
Indeed, \(C(k+1)/C(k)\ge C(n)/C(n-1)\) throughout this range,
with equality at \(k=n-1\). For \(k<0\) or \(k\ge n\) the
corresponding distribution-function bounds follow from the support. At a common quantile we obtain
\(K_{n-1}\le K_n\le K_{n-1}+1\). Given these weights, couple the
lengths as above and draw a word uniformly from each resulting
weight and length class. Each word then has the required marginal
probability. Both \(J_v\) and the unused weight \(v-K_v\) increase
weakly along this coupling. Since \(s_\ell\) is nonincreasing in
\(\ell\), the event \(v-K_v\ge s_{J_v}\) persists from one budget to
the next. Its probability is exactly \eqref{eq:threshold}.

To close the induction on the ratios, use
\((1-z)S_t=z-z^{t+1}\) and \(P_{\ell+1}=S_{b_\ell}P_\ell\)
for \(\ell<H\) to write
\[
 (1-2z)F=1-z-z\sum_{\ell=0}^H z^{b_\ell}P_\ell.
\]
The terminal term \(b_H=0\) arises from the last layer. Dividing
by \(1-z\) and comparing coefficients at \(n+1\) gives
\begin{equation}\label{eq:rec}
 \frac{C(n+1)}{C(n)}
 =2-\frac{\sum_{\ell=0}^H C_\ell(n-b_\ell)}{C(n)}
 \qquad(n\ge0).
\end{equation}
The empty word ensures \(C(n)>0\). The induction starts with the
single ratio \(C(1)/C(0)\). Given the ordering through \(n\), take
\(s_\ell=b_\ell\) in the coupling for \(n-1,n\). The fraction subtracted
in \eqref{eq:rec} increases weakly, so
\(C(n+1)/C(n)\le C(n)/C(n-1)\). This also constructs the coupling at
every pair of adjacent budgets. For the empty host \(P_0=1\) and
\(b_0=0\), and the recurrence gives \(C(n)=1\).
\end{proof}
\section{Positive rank differences and cumulative counts}\label{sec:arms}

For an increasing arm, we need to compare two growth rates. If
\(d_{N-1}>0\) and \(N\ge3\), the ratio \(d_N/d_{N-1}\) is at most
\(G_j/G_{j-1}\) at every earlier rank \(1\le j\le N-2\).
To prove this inequality we compare both ratios with the cumulative
counts of a shorter arm.

Normalize a noninterval increasing arm to
\(X=1^\alpha(t_1,\ldots,t_h)\), with \(h\ge1\). Write
\(G_n=G_n(X)\) and \(d_n=d_n(X)\).
Take as reference the decreasing host \((t_h,\ldots,t_2)\),
with its \(P_\ell,C_\ell,C\) as in Section~\ref{sec:transport}.
Set \(P_h=\prod_{i=1}^h S_{t_i}\) for the one additional layer.
For \(h=1\) the reference host is empty.

\begin{theorem}[Comparison with a shorter arm]\label{thm:sandwich}
For \(m\ge0\), \(d_{m+2}/C(m)\) is nonincreasing, and \(G_m/C(m)\)
is nondecreasing.
\end{theorem}
\begin{proof}
For a series \(F\) with initial coefficients \(1,1\), write
\((1-z)F=1+z^2 L(z)\). Under an append \(F\mapsto 1+S_tF\)
its signed difference series \(L\) becomes
\([t-1]+S_t L\). Indeed
\(1-z+S_t=1+z^2[t-1]\). The two possible initial series
\((1-z)^{-1}\) and \([\alpha+1]\) give initial \(L=0\) and
\(L=-z^{\alpha-1}\), respectively. Starting from the infinite series gives
\[
 L_\infty(z)=\sum_{\ell=0}^{h-1}[t_{h-\ell}-1]P_\ell(z),\qquad
 \coef m{L_\infty}=C(m)-\sum_{\ell=0}^{h-1}C_\ell(m-t_{h-\ell}+1).
\]
Here \(P_\ell\) consists of the last \(\ell\) append factors, with
\(\ell\le h-1\) as required for the reference host. In the finite arm,
the initial difference series contributes \(-T\), where
\(T=z^{\alpha-1}P_h\). Thus \(d_{m+2}=((L_\infty)_m-T_m)_+\).

The thresholds \(t_{h-\ell}-1\) decrease as \(\ell\) increases,
so \eqref{eq:threshold} implies that \((L_\infty)_m/C(m)\) is nonincreasing.

For \(i=h-\ell\), write the corresponding summand
\(L_i=P_\ell[t_i-1]\) of \(L_\infty\) and split
\[
 T=P_\ell z^{\alpha+i-1}\prod_{j=1}^i[t_j].
\]
The second factor on the right has nondecreasing coefficients on
\(0,\ldots,t_i-2\), the support of \([t_i-1]\). In that range its
difference is an entering coefficient of the convolution with
\([t_i]\), since the departing index is negative. The coefficients of
\([t_i-1]\) and of this second factor are therefore in ratio order,
including the endpoints. Convolution with \(P_\ell\) preserves the
order, so \((L_i)_mT_{m'}\ge(L_i)_{m'}T_m\) for \(m<m'\). Summation
over \(i\) gives the same inequality for \(L_\infty,T\).

The positive support of \(L_\infty\) is an initial interval because
the summand indexed by \(\ell\) starts at \(\ell\) and extends
through at least that index. On the support,
\((1-T_m/(L_\infty)_m)_+\) is nonincreasing. Multiplying it by
\((L_\infty)_m/C(m)\) and adjoining zeros proves the first assertion.

Iterating \eqref{eq:append} also gives
\[
  G_m=C(m)+\coef m{\frac{P_h[\alpha+1]}{1-z}}.
\]
Let \(V_m\) be the added coefficient. For each \(\ell<h\), \(V\)
is the convolution of \(C_\ell\) with the nonnegative coefficients of
\([\alpha+1]\prod_{i=1}^{h-\ell}S_{t_i}\).
The sequence \(C_\ell\) is log-concave by convolution with the
all-ones sequence. Its translate comparison gives
\(V_{m+1}C_\ell(m)\ge V_m C_\ell(m+1)\), including the zero indices.
Summing this inequality proves that \(V_m/C(m)\) is nondecreasing,
since \(C(m)=\sum_\ell C_\ell(m)\).
\end{proof}

\begin{corollary}[The single-arm comparison]\label{cor:arm}
For every increasing arm, \(N\ge3\), and \(1\le j\le N-2\),
\begin{equation}\label{eq:arm}
 G_j d_{N-1}\ge G_{j-1}d_N .
\end{equation}
Also \(d_N\le2d_{N-1}\) in this range of \(N\).
\end{corollary}
\begin{proof}
An interval arm has \(d_N=0\) for \(N\ge2\). For any other arm,
Theorem~\ref{thm:sandwich} gives \(d_N=0\) whenever \(d_{N-1}=0\).
We may therefore assume that both differences are positive. The theorem
and log-concavity of \(C\) give
\[
 \frac{d_N}{d_{N-1}}
 \le \frac{C(N-2)}{C(N-3)}
 \le \frac{C(j)}{C(j-1)}
 \le \frac{G_j}{G_{j-1}}.
\]
The first \(C\)-ratio is at most two by \eqref{eq:rec}.
\end{proof}

We shall use \eqref{eq:arm} with \(N=n+1\) and \(j=n+1-r\) for
\(n\ge2\), \(r\ge3\). The convention \(G_j=0\) for \(j<0\) gives
the inequality also at \(j\le0\).
\section{Weighted positive differences of products}\label{sec:product}

For increasing arms \(X,Y\), put
\[
 \begin{gathered}
 A=F_X,\quad B=F_Y,\quad a=\wt(X),\quad b=\wt(Y),\\
 P=AB,\quad p=q_Xq_Y,\quad
 u_n=(p_n-p_{n-1})_+,\quad U_n=2^n u_n\quad(n\ge0).
 \end{gathered}
\]
Thus \(p_n=P_{a+b-n}\) and \(u_0=p_0=1\). We write
\(q_{X,i}=\coef{i}{q_X}\) and use the operator
\((\mathcal D v)_i=(\Delta v_i)_+\) on coefficient sequences.
Proposition~\ref{prop:count-induction} asserts that \(U\) is unimodal.
At the peak of size \(M\) we will compare the gain \(p_{M+k}-p_k\) with
\(u_k\) in \eqref{eq:joining}. Once \(k\ge\firstmode(U)\), unimodality
gives \(u_{k+i}\le2^{-i}u_k\) for \(i\ge1\), and hence
\[
 p_{M+k}-p_k\le\sum_{i=1}^M u_{k+i}\le u_k .
\]

Every noninterval arm has at least two positions and a part of size
at least two, so it contains both \((2)\) and \((1,1)\).
With two such arms, \(P\) begins with \(1,2,5\) and fails log-concavity
at rank one. Chen (\href{https://doi.org/10.2139/ssrn.7385218}{SSRN~7385218})
proved that its coefficients are log-concave
from their first strict descent onward. To study the positive
differences of the reversed product, take \(X=Zt\) to be the heavier
arm. Then \(p=[t]q_Zq_Y+z^a q_Y\) by \eqref{eq:append}. The added
term begins after the first mode of \(p\), so it suffices to prove
that the interval convolution \([t]q_Zq_Y\) has unimodal weighted
positive differences.

\subsection{Coefficients and positive differences}

\begin{lemma}[Ordinary unimodality and the first mode]\label{lem:ordinary}
The products \(P,p\) are unimodal. If neither arm polynomial is an
interval, then \(\firstmode(p)<\max(a,b)\).
The sequence \(u\) is unimodal and has initial positive support.
\end{lemma}
\begin{proof}
By Lemma~\ref{lem:layers}, the height intervals give
\(P_0\le\cdots\le P_{\lceil(a+b)/2\rceil}\). If \(a=b\) and neither
arm is an interval, we can go one rank farther. From \(a\) to \(a+1\)
the product of the two lowest intervals loses one, while the product
of the two second intervals \([2,a-1]\) gains one. Every other pair of
layers has endpoint sum at least \(2a+1\) and contributes a nonnegative
difference. Thus \(P_{a+1}\ge P_a\).

With an interval factor, \(P\) is unimodal by Lemma~\ref{lem:layers}.
Otherwise take \(a\ge b\) and write \(X=Zt\). In \(P=B+S_t F_ZB\) the
second summand is unimodal by induction on the number of parts and
convolution with an interval. The preceding bounds give nondecrease
through \(b+1\), and from that rank onward \(P\) agrees with the second
summand. Its last mode is therefore greater than \(b\), which puts the
first mode of \(p\) below \(a\).

Let \(v\) be nonnegative and unimodal, with last mode \(m\), and
suppose \(\mathcal Dv\) is unimodal. Through \(m\), the coefficients of
\(\mathcal D([t]v)\) agree with \([t]\mathcal Dv\), a unimodal sequence
by Lemma~\ref{lem:layers}. At \(m+1\) the new coefficient is at most
the corresponding coefficient of \([t]\mathcal Dv\), so a fall already
present in the prefix cannot be followed by a rise across this join.
Beyond \(m\), while \(v_n-v_{n-t}>0\), the entering coefficient
decreases weakly and the departing coefficient increases weakly.
Once \(n-t\) reaches a mode, the difference is nonpositive. Thus
\(\mathcal D([t]v)\) is unimodal.

Now delete the last part of a heavier arm. If neither
factor is an interval, orient \(a\ge b\) and write
\(p=[t]q_Zq_Y+z^a q_Y\). Up to the first mode of \(p\) the
last term is absent. Its positive differences are therefore a
prefix of the window's positive differences, followed by zero.
If one factor is an interval, it suffices to treat \(q\) for a
single arm and then apply the observation. A single noninterval
arm has first reversed mode below its weight, since its last
coefficient is one and its maximum exceeds one. Apply
\eqref{eq:append} and the same truncation argument to that arm.
A single interval starts the induction. Since \(u_0=1\),
unimodality gives initial positive support.
\end{proof}

\begin{lemma}[Rising log-concavity]\label{lem:rising-lc}
For a reversed arm or a reversed two-arm product \(v\),
\(\mathcal D^i v\) is unimodal and log-concave through its last mode
for every \(i\ge0\), with \(\mathcal D^0 v=v\).
\end{lemma}
\begin{proof}
Let \(v\) be nonnegative and unimodal, log-concave through its last
mode \(m\), and put \(w=[t]v\). For \(n+1\le m\), truncating \(v\) at
\(m\) gives a log-concave sequence whose convolution agrees with \(w\)
through \(n+1\). For \(n\ge m\), a positive difference
\(\Delta w_{n+1}=v_{n+1}-v_{n+1-t}\) places the departing index before
the first mode of \(v\), while the entering index is past the last.
Hence \(0<\Delta w_{n+1}\le\Delta w_n\), and
\[
 w_n^2-w_{n-1}w_{n+1}
 =w_n(\Delta w_n-\Delta w_{n+1})
   +\Delta w_n\,\Delta w_{n+1}\ge0 .
\]
If \(\Delta w_{n+1}\le0\) at a rank on the rising side, including the
last mode, then \(w_n\ge w_{n-1},w_{n+1}\) gives the same inequality.

Next suppose the support of \(\mathcal Dv\) is \([0,L]\).
Through \(L\), \(\mathcal D([t]v)\) agrees with \([t]\mathcal Dv\).
At every subsequent step, the entering difference is nonpositive and
the departing difference is nonnegative until the window has passed
\([0,L]\), after which the positive part is zero. Hence the tail
starting at \(L\) is nonincreasing. The order-zero window argument
applied to \(\mathcal Dv\) gives log-concavity on the rising prefix
of \(\mathcal D([t]v)\). The joining inequality, when needed at index
\(L\), follows from the fact that its two neighbors are at most its
value there.

It follows that \(\mathcal D^2([t]v)\) equals
\(\mathcal D([t]\mathcal Dv)\) through \(L\), followed by zero.
For higher iterates the same identity holds with both operator powers
increased by one at each step, since taking positive differences
commutes with truncating and filling the tail with zeros. Truncation
preserves unimodality and log-concavity through the last mode. Induction
on the order, applying the window assertion to \(\mathcal Dv\) at the
preceding order, proves the result for each iterate of the convolution.
The initial-support condition holds for the inputs at each step because
their positive differences are unimodal and have value \(v_0>0\) at zero.

For a single reversed arm, an interval starts the induction.
In a noninterval append the terminal unit \(z^a\) changes
the final difference of the window from \(-1\) to zero and
does not change positive differences. All inequalities
before the first mode transfer from the window, and those
on the maximal plateau are automatic. For two noninterval
arms use the heavier deletion in Lemma~\ref{lem:ordinary}.
Up to the first mode the new product agrees with the old
window, and the positive differences are the corresponding
truncation. The same reasoning handles order zero through the
first mode and then its maximal plateau. An interval factor
is exactly a uniform convolution.
\end{proof}

\subsection{Source polynomials and the base cases}

For an increasing arm \(X\) of weight \(a\), let \(\kappa(X)\) be
the least weight of a composition outside its ideal. The one-part
composition \((\max(X)+1)\) shows that \(\kappa(X)\le\max(X)+1\).
Define the
forward source \(\mathcal E_X\) and its reversal \(R_X\) by
\begin{equation}\label{eq:sources}
 \mathcal E_X=1-z+(2z-1)F_X,\qquad
 R_X=z^{a+1}\mathcal E_X(1/z)=(2-z)q_X-z^a+z^{a+1}.
\end{equation}
At rank \(i\ge2\), the coefficient of the forward source
\(\mathcal E_X\) is \(2a_{i-1}(X)-a_i(X)\), the slack in
Lemma~\ref{lem:count-bounds}. Its reversal \(R_X\) has the same
coefficients as \((2-z)q_X\) below degree \(a\).

\begin{lemma}[The first missing weight]\label{lem:source}
The polynomials \(\mathcal E_X\) and \(R_X\) are nonnegative and
unimodal. The support of \(\mathcal E_X\) is \([\kappa(X),a+1]\).
For \(X=1^\alpha(t_1,\ldots,t_h)\) with \(\alpha,h\ge1\) and
\(t_i\ge2\),
\begin{equation}\label{eq:kappa}
 \kappa(X)=\min\bigl(\{h+2\}\cup\{t_i+h-i+1:1\le i\le h\}\bigr).
\end{equation}
\end{lemma}
\begin{proof}
The empty arm has source \(z\) and \(\kappa(\varnothing)=1\).
For an append to \(X\) of weight \(a\), the recurrence is
\(\mathcal E_{Xt}=S_t\mathcal E_X+z^{t+1}\), hence
\(R_{Xt}=[t]R_X+z^a\). This proves nonnegativity by induction.

For a nonempty arm, the source coefficients at ranks zero and one
vanish. There are \(2^{n-1}\) compositions of weight \(n\ge1\).
Below the first missing weight the rank counts double at each step
from rank one onward. At \(n=\kappa(X)\) we have
\(a_n(X)<2a_{n-1}(X)\), so the source starts at \(k=\kappa(X)\).
Since \(k\le\max(X)+1\le t+1\), the window and the new term at \(t+1\)
have interval support together, starting at \(\min(k+1,t+1)\) and
ending at \(a+t+1\). This proves the support claim by induction.
Starting with \(k=2\) for a nonempty string of ones and iterating
the starting-index rule gives \eqref{eq:kappa}.

For unimodality, the window \(S_t\mathcal E_X\) is unimodal by
induction. Its differences through \(t+1\) are nonnegative because
the departing source coefficient has index at most zero. Adding
at \(t+1\) increases the difference there and decreases the one
at \(t+2\), which was \((\mathcal E_X)_{t+1}-(\mathcal E_X)_1\).
For a nonempty prefix \((\mathcal E_X)_1=0\). If
\((\mathcal E_X)_{t+1}>0\), the change by \(-1\) preserves the
nonnegative edge. Otherwise the support bound
\(k\le t+1\) implies \(t+1>a+1\), so all subsequent window
differences are nonpositive. An empty prefix gives
\(zS_t+z^{t+1}\), which is unimodal directly.
\end{proof}

\begin{lemma}[Cumulative ratios for the reversed source]\label{lem:source-ratio}
For every increasing arm, the cumulative sequence
\(C_n=\sum_{i\le n}(R_X)_i\) is log-concave. Consequently, for
\(w=[t]R_X\) with \(t\ge2\),
\begin{equation}\label{eq:source-ratio}
 w_{k+1}w_j\ge w_kw_{j+1}\qquad(0\le k\le t-2,\ j\ge k).
\end{equation}
\end{lemma}
\begin{proof}
In the append \(R_{Xt}=[t]R_X+z^a\), set \(v_n=([t]C)_n\),
\(w_n=v_n-v_{n-1}\). By the log-concave convolution property of
Section~\ref{sec:transport}, \(v\) is log-concave.
The new cumulative sequence is \(v_n+\mathbf1_{n\ge a}\).
For an empty prefix it is \(2,3,\ldots,t+1\) with a constant
tail. For \(a\ge t\) the source support gives
\(w_{a-1}-w_a=(R_X)_{a-t}\ge1\), since
\(\kappa(X)\le t+1\) and \((R_X)_a=0\). By unimodality of
\(w\), it is nonincreasing from \(a-1\) onward. At the center
\(a-1\) the new log-concavity determinant is
\(v_{a-2}(w_{a-1}-w_a-1)+w_{a-1}^2\ge0\).
If \(1\le a<t\), both \(w_{a-1}\) and \(w_a\) equal the source
mass \(m=F_X(1)\ge a+1\), and \(w\) is nonincreasing from \(a-1\).
The determinant at \(a-1\) is then
\(m^2-v_{a-2}\ge m^2-(a-1)m\ge0\).
At \(a\) the added step increases the determinant by
\(2v_a+1-v_{a-1}\ge0\), while at the later centers it
increases by \(2v_n-v_{n-1}-v_{n+1}\ge0\).
At earlier centers the three coefficients are unchanged.

To prove \eqref{eq:source-ratio}, put \(r=C_{k+1}/C_k\),
since \(w_i=C_i\) at \(i=k,k+1\). For \(k\le j\le k+t\),
subtract \(C_{j+1-t}\ge r C_{j-t}\) from
\(C_{j+1}\le rC_j\), using zero extension at the lower boundary.
For \(i\ge k+1\), log-concavity gives
\((R_X)_{i+1}C_{i-1}\le (R_X)_i C_i\) and
\(C_i/C_{i-1}\le r\).
Summing \((R_X)_{i+1}\le r(R_X)_i\) over the
appropriate window proves the result for \(j\ge k+t\).
\end{proof}

\begin{lemma}[Base cases for the weighted differences]\label{lem:weighted-bases}
The sequence \(U\) is unimodal if one arm is an interval arm, or if
an arm of maximum weight has only one part greater than one after
normalization.
\end{lemma}
\begin{proof}
Suppose \(q_X=[\alpha+1]\). For \(\alpha\ge b\), the first
product mode is \(b\). The positive differences through \(b\)
are \(q_{Y,j}\). Their weighted values increase by
Lemma~\ref{lem:count-bounds}. For \(\alpha<b\), the difference
at \(b\) is \(1-q_{Y,b-\alpha-1}\le0\), and later
differences are nonpositive. Thus the first product mode is
below \(b\). Up to it the polynomial \((2-z)p\) agrees with
\([\alpha+1]R_Y\), which is unimodal. Up to that mode its differences have the same signs as those of \(U\), since
\[
 [(2-z)p]_j-[(2-z)p]_{j-1}=2u_j-u_{j-1}
 \quad(1\le j\le\firstmode(p)).
\]
After the first mode \(U\) vanishes. This also treats an empty arm.

We may now suppose \(X=1^\alpha t\)
with \(t\ge2\) and \(a=\alpha+t\ge b\), and \(Y\) noninterval.
The arm \(1^{t-1}(\alpha+1)\) has the same rank polynomial
\(1+z[\alpha+1][t]\), with the two interval factors interchanged.
We may therefore take \(\alpha\ge t-1\). Put \(V=[t]q_Y\) and
\(W=[t]R_Y\).
By Lemma~\ref{lem:ordinary}, the first mode \(\sigma\) of \(p\)
is below \(a\), and up to \(\sigma\) the positive differences
are
\(u_j=V_j-V_{j-\alpha-1}\).
In particular they are \(V_j\) through \(\alpha\), within this
range. The coefficients of \([t]B\) grow by a factor at most
two in forward order. For \(n<t\), use
\(B_n\le\sum_{i<n}B_i\). For \(n\ge t\), sum the
individual bounds \(B_i\le2B_{i-1}\) over the window.
Both bounds are in Lemma~\ref{lem:count-bounds}.
It follows that \(2^j V_j\) is nondecreasing up to the end of
its support.

For \(\alpha+1\le j\le\sigma\), set \(k=j-\alpha-1\).
Then \(k\le t-2\), and \(\sigma\le b+t-1\) because all
window differences beyond \(\deg V\) are nonpositive before
the terminal correction. Consequently \(k<b\).
Using \eqref{eq:sources} gives
\[
 (2-z)V=W+z^b-z^{b+t},\qquad
 2u_j-u_{j-1}=W_j-W_k+\mathbf1_{j=b}.
\]
The coefficients \(W_k\) are positive since \(k<t\)
and the source starts at rank zero. By
\eqref{eq:source-ratio}, \(W_{k+\alpha+1}/W_k\) cannot rise
over these indices. The differences \(W_j-W_k\) cannot
go from negative to positive. If one has already become
negative before \(j=b\), it is at most \(-1\) at \(b\), so
adding one still leaves a nonpositive difference. Thus the weighted
sequence can change direction only once before it reaches zero.
\end{proof}

\subsection{A local inequality for interval convolution}

In \(q_{Xt}=[t]q_X+z^{a+t}\), the interval factor is the term that
can move the first mode. To follow its effect on \(U\), we take
differences of the coefficients of \((2-z)p\).

For any polynomial
\(v\) write
\[
 E_v(n)=2v_n-v_{n-1},\qquad T_v(n)=E_v(n+1)-E_v(n).
\]
For \(p=q_Xq_Y\), abbreviate \(E_n=E_p(n)\). At a fixed reverse
index \(n\), write \(T=T_p(n)\) and let \(j=a+b-n\) be its
forward index. Then \(T=P_{j+1}-3P_j+2P_{j-1}\).

\begin{lemma}[Two increasing factors]\label{lem:tilt}
For two noninterval arms, the sequence \(2^np_n\) is unimodal.
At any \(n\) for which \(E_n\ge0\) and both forward factors increase
strictly from \(j=a+b-n\) to \(j+1\), one has \(E_{n+1}\le E_n\).
\end{lemma}
\begin{proof}
For a noninterval arm, the coefficients of \(\mathcal E_Y\)
are zero through rank two. At every later rank where \(B\)
is still strictly rising, the source increment is
\(2d_{i-1}(Y)-d_i(Y)\ge0\) by Corollary~\ref{cor:arm}.
A positive \(d_i(Y)\) in this range has a positive predecessor,
except at the initial rank two, by Theorem~\ref{thm:sandwich}.
The same assertion holds for \(A\).

Define \(e_j=2P_{j-1}-P_j\) and \(L_j=(\mathcal E_Y A)_j\).
Then \(e_j=L_j-\Delta A_j\). At consecutive ranks
\(j,j+1\ge2\) of common strict increase,
\[
 L_i=\sum_{s=3}^i(\Delta\mathcal E_Y)_s G_{i-s}(X).
\]
Corollary~\ref{cor:arm}, applied with \(N=j+1\) and
with index \(j+1-s\), gives
\(L_{j+1}d_j(X)\ge L_jd_{j+1}(X)\).
The term at \(s=j+1\) is nonnegative, and \(L_{j+1}\ge L_j\).
If \(e_{j+1}<0\) the ratio comparison gives \(e_j<0\).
If \(e_{j+1}\ge0\), then \(L_{j+1}>0\) and
\[
 L_{j+1}(e_{j+1}-e_j)
 =L_{j+1}d_j-L_jd_{j+1}+(L_{j+1}-L_j)e_{j+1}\ge0 .
\]
If either factor ceases to rise at some rank \(j\ge3\), orient it as
\(A\). Its positive forward differences from rank two onward have
initial support by Theorem~\ref{thm:sandwich}, so a strict rise cannot
resume.
The identity for \(e_j\) gives \(e_j\ge0\) at every such
rank. The initial values \(P_0,P_1,P_2\) are \(1,2,5\),
so \(e_1=0,e_2=-1\). These facts give unimodality of the
coefficients of \(P(z/2)\), and hence of their scaled reversal
\(2^np_n\). Since \(E_n=e_{j+1}\),
the inequality \(e_j\le e_{j+1}\) proves the other conclusion.
The possible edge \(j=1\) is excluded by \(E_n\ge0\) and \(e_2=-1\).
\end{proof}

In the induction of Proposition~\ref{prop:count-induction}, we apply
the following comparison to the product \(p\) for the shortened pair
of arms, whose weighted positive differences are already unimodal.
For the window \([t]p\), the difference at rank \(n\) is \(p_n-p_k\),
where \(k=n-t\). The comparison controls the weighted differences
once the entering coefficient is past the last mode of \(p\).

\Needspace{9\baselineskip}
\begin{condition}[Local comparison for two noninterval arms]\label{cond:local}
Suppose \(U\) is unimodal for \(p=q_Xq_Y\). Put \(m=\lastmode(p)\)
and \(\ell=\firstmode(U)\). The local comparison holds if, for all
\(n>m\) and \(0\le k<\min(n,\ell)\), the inequalities
\begin{equation}\label{eq:local-hypotheses}
 p_n>p_k,\qquad p_{n+1}>p_{k+1},\qquad 0\le E_n\le E_k
\end{equation}
imply
\begin{equation}\label{eq:local}
 E_{n+1}-E_n\le E_{k+1}-E_k .
\end{equation}
\end{condition}

The positive differences of \(p\) vanish after its first mode, so
\(k<\ell\) puts both \(k\) and \(k+1\) at or before that mode.
The right side of \eqref{eq:local} is then
\(2u_{k+1}-u_k\ge0\) by the increase of \(U\) through \(\ell\).
When both forward factors increase from \(j\) to \(j+1\),
Lemma~\ref{lem:tilt} makes the left side of \eqref{eq:local}
nonpositive.

\begin{lemma}[Preservation under interval convolution]\label{lem:window}
For a finite nonnegative sequence \(p\) with \(p_0>0\), suppose that
\(p\), \((2^np_n)_n\), and \(U\) are unimodal. If
\eqref{eq:local-hypotheses} implies \eqref{eq:local} for every
\(n>\lastmode(p)\) and \(0\le k<\min(n,\firstmode(U))\), then the
weighted positive differences of
\([t]p\) are unimodal for every \(t\ge1\).
\end{lemma}
\begin{proof}
Write \(w_n=p_n-p_{n-t}\) and \(f_n=2^n(w_n)_+\).
The positive support of \(u\) is an initial interval since \(U\) is
unimodal and \(u_0>0\). Hence \(p\) rises strictly to its first mode.
The indices with \(w_n>0\) likewise form an initial interval. Up to the
last mode of \(p\), the first equality \(p_n=p_{n-t}\) places both
indices on the modal plateau. After the last mode, positivity requires
\(n-t\) to precede the first mode, and \(p_n-p_{n-t}\) decreases weakly
until it reaches zero.

For \(n\le m=\lastmode(p)\), \(f_n\) agrees with the convolution
of \(U\) with \((1,2,\ldots,2^{t-1})\). It is unimodal. Indeed, the positive and negative parts of
\(\Delta U\) have ordered supports and hence are in ratio order.
Convolving both parts with \((1,2,\ldots,2^{t-1})\) preserves that order
by the Toeplitz-minor argument of Section~\ref{sec:transport}. Their
difference can therefore change sign only from positive to negative.
At \(m+1\), \(f_{m+1}\) is at most the corresponding coefficient of
this convolution. There is also the identity \(f_n=2^np_n\) for
\(n<t\), with \(f_t\le2^tp_t\).

Choose the first comparison if \(m\ge t-1\) and the second
otherwise. With \(c=\max(m,t-1)\), an earlier strict fall in
\(f\) through \(c\) therefore gives \(f_{c+1}\le f_c\).

For the remaining edges, suppose \(n\ge c+1\),
\(f_n\le f_{n-1}\), and \(w_n,w_{n+1}>0\).
Set \(k=n-t\ge0\). In particular \(n>m\) and the outgoing
index \(k\) precedes the first mode. If \(k\ge\ell\),
successive halvings after the first mode of \(U\) give
\(w_n\le\sum_{i>k}u_i\le2u_{k+1}\).
Since \(p_{n+1}\le p_n\),
\[
 w_{n+1}=w_n-u_{k+1}-(p_n-p_{n+1})\le w_n/2 .
\]
If \(k<\ell\), the supposed fall is
\(E_n-E_k=2w_n-w_{n-1}\le0\). Also \(E_{k+1}>0\) and
\(E_{k+1}-E_k=2u_{k+1}-u_k\ge0\).
When \(E_n<0\), unimodality of \(2^np_n\) implies
\(E_{n+1}\le0\). When \(E_n\ge0\), use \eqref{eq:local}.
Both give \(2w_{n+1}-w_n\le0\).
It follows that once \(f\) decreases, it remains nonincreasing to the
end of its positive support and is zero thereafter.
\end{proof}

\subsection{Bounds from the parts of an arm}

The least weight of a word outside an arm ideal places a lower bound on
the initial interval where the reversed coefficients rise strictly.

\begin{lemma}[An initial strict range]\label{lem:head-range}
For a normalized nonempty arm \(Y=1^\beta(v_1,\ldots,v_g)\), define
\[
 M_Y=\beta+\sum_{i=1}^g(v_i-1),\quad
 B_Y=\max\{\beta,v_1-1,\ldots,v_g-1\},\quad
 s_Y=\min(\lfloor M_Y/2\rfloor,M_Y-B_Y).
\]
Then \(\Delta q_{Y,i}\ge1\) for \(0\le i\le s_Y\).
If \(K=\kappa(Y)\ge3\), then
\[
 s_Y\ge S_K:=\left\lfloor\frac{1+(K-2)(K-1)/2}{2}\right\rfloor,\qquad
 \wt(Y)\ge K(K-1)/2 .
\]
Moreover, with \(J_Y=[\beta+1]\prod_i[v_i]\),
\(\Delta q_{Y,i}\le (J_Y)_i\) at every rank.
\end{lemma}
\begin{proof}
Iterating \eqref{eq:append} gives
\(q_Y=J_Y+\sum_{i=1}^g z^{w_i}\prod_{j>i}[v_j]\), where
\(w_i=\beta+\sum_{j\le i}v_j\).
Apart from \(J_Y\), every summand has endpoint sum at least \(M_Y\),
so contributes nonnegative differences through
\(\lfloor M_Y/2\rfloor\).
In an interval product of degree \(M\) with largest factor
degree \(B\), the differences are strictly positive through
\(\min(\lfloor M/2\rfloor,M-B)\). For the induction, remove a largest
factor. The remaining product is \(Q\) of degree \(L=M-B\),
with largest factor degree \(B'\le B\).
The empty product starts the induction. The difference at \(i\) is
\(Q_i-Q_{i-B-1}\).
For \(i\le\min(\lfloor M/2\rfloor,L)\), the first term
is positive. If \(j=i-B-1\ge0\), then \(i+j\le L-1\)
and \(j<\min(\lfloor L/2\rfloor,L-B')\).
Symmetry and the strict increasing range for \(Q\)
give \(Q_i>Q_j\).

By \eqref{eq:kappa}, \(g\ge K-2\) and the last \(K-2\)
parts dominate \(2,\ldots,K-1\). This yields the weight
bound and \(M_Y\ge1+(K-2)(K-1)/2\).
After removing any one interval factor, the total degree
left is at least \(1+(K-3)(K-2)/2\), which is at least
\(S_K\). Hence \(s_Y\ge S_K\).
To bound \(\Delta q_{Y,i}\) from above, start with \([\beta+1]\).
At each append, convolution with the nonnegative interval preserves
the bound up to the rank preceding the terminal correction. At the
new arm weight the last two coefficients of its reversed polynomial
are one, so the difference there is zero and the next is negative.
This yields \(\Delta q_{Y,i}\le(J_Y)_i\) at every rank.
\end{proof}

\begin{lemma}[A bound from the largest part]\label{lem:capacity}
Let \(X=Zt\) be noninterval, suppose \(A_{j+1}\le A_j\), and
put \(r=a-j\) and \(C_0=\Delta q_{X,r+1}-\Delta q_{X,r}\).
If \(C_0>0\), then
\[
 C_0\le q_{X,r+\delta}\qquad(1\le\delta\le t).
\]
\end{lemma}
\begin{proof}
For \(r\le-2\) the curvature is zero. For \(r=-1\),
\(C_0=1\) and \(q_{X,\delta-1}>0\).
Otherwise the last two coefficients of \(q_X\) are one,
so \(C_0>0\) and \(\Delta q_{X,r}\ge0\) force \(r\le a-2\).
Put \(Q=q_Z\). The correction at \(a\) is absent through
\(r+1\), and
\[
 \Delta q_{X,r}=Q_r-Q_{r-t},\qquad
 C_0=\Delta Q_{r+1}-\Delta Q_{r+1-t}.
\]
If \(\Delta Q_{r+1}\le0\), then \(\Delta Q_{r+1-t}<0\).
Unimodality of \(Q\) would give \(Q_r<Q_{r-t}\), a
contradiction. Thus \(\Delta Q_{r+1}>0\), and the earlier
difference at \(r+1-t\) is nonnegative. Consequently
\(C_0\le Q_{r+1}\). For every \(1\le\delta\le t\), the coefficient \(Q_{r+1}\) occurs in
\(([t]Q)_{r+\delta}\), which is at most \(q_{X,r+\delta}\).
\end{proof}
\subsection{When both factors are nonincreasing}

Since \(E_{k+1}-E_k\ge0\) for \(k<\firstmode(U)\), it suffices here
to show \(T_p(n)\le0\). We begin with the height intervals of the
reversed source.

\begin{lemma}[Reflected source intervals]\label{lem:source-reflection}
For any nonempty increasing arm of weight \(a\), the height intervals
of \(R_X\) have endpoint sum at most \(a-1\).
\end{lemma}
\begin{proof}
In forward order the claim is that the endpoint sums of
\(\mathcal E_X\) are at least \(a+3\). This holds for a
single part of size \(a\), where
\(\mathcal E_X=z^2[a-1]+2z^{a+1}\) (the first term is absent
for \(a=1\)). For an append, write \(c\) for the old weight
and \(W=S_t\mathcal E_Z\). The old interval bounds give
endpoint sums at least \(c+t+4\) for every interval product
in \(W\). Only the unit added at \(t+1\) can violate the coefficient
comparison about \((c+t+3)/2\). It can do so only if \(c\ge t\),
and its reflected index is then \(c+2\). Write \(e_i=(\mathcal E_Z)_i\).
Reflecting the interval \([1,t]\) about \((c+3)/2\), and
pairing terms where the two intervals overlap, gives
\(\sum_{i=1}^t e_i\le\sum_{i=c-t+3}^{c+2}e_i\). Hence
\[
 W_{c+2}-W_{t+1}
 =\sum_{i=c-t+2}^{c+1}e_i-\sum_{i=1}^t e_i
 \ge e_{c-t+2}.
\]
If \(e_{c-t+2}=0\), the first sum includes the entire positive
source support, and the second omits \(e_{c+1}=2\).
In either case \(W_{c+2}\ge W_{t+1}+1\). The reflected inequality
therefore holds after the append. Since the source is unimodal
(Lemma~\ref{lem:source}), its height intervals have the claimed
endpoint sums.
\end{proof}

\begin{lemma}[The first half with two nonincreasing factors]\label{lem:both-nonrise}
Let both arms be noninterval and let the two forward factors be
nonincreasing at \(j+1\). If \(2j\le a+b+2\), then
\(T_p(a+b-j)\le0\).
\end{lemma}
\begin{proof}
Write \(C_0,r\) as in Lemma~\ref{lem:capacity}, and
\(K=\kappa(Y)\). By \eqref{eq:sources},
\[
 (2-z)p=R_Yq_X+z^b(1-z)q_X .
\]
All components of \(R_Yq_X\) in the interval expansion
have endpoint sum at most \(a+b-1\), by
Lemmas~\ref{lem:layers} and \ref{lem:source-reflection}.
Put \(n=a+b-j\). Each component contributes a nonnegative drop
from its coefficient at \(n\) to that at \(n+1\).
For the lowest source interval \([0,b+1-K]\) convolved with
\(q_X\), this drop is
\(q_{X,r+K-1}-q_{X,n+1}\).
In particular, for \(j\le b\), \(q_{X,n+1}=0\) and
\(T_p(n)\le C_0-q_{X,r+K-1}\).

For \(j\le\min(a,b)\), orient \(X\) with the larger
last part \(t\), so \(K-1\le t\), and apply
Lemma~\ref{lem:capacity} whenever \(C_0>0\).
For the other cases orient \(a\le b\).
If \(a+2\le j\le b\), then \(C_0=0\). If \(j=a+1\le b\),
then \(C_0=1\), and the nonincrease of \(B\) at \(j+1\)
implies \(K\le j+1\), since the full composition counts
strictly increase there. Thus \(q_{X,K-2}\ge1\).
Finally, if \(j>b\), the half-range inequality forces
\(a=b\) and \(j=a+1\). Again \(C_0=1\).
The lowest source drop is \(q_{X,K-2}-1\).
If the normalized \(Y\) has \(g\) parts greater than one,
then \(3\le K\le g+2\) by \eqref{eq:kappa} and \(2g\le a-1\).
Thus
\(1\le K-2\le a-2\).
All coefficients of \(q_X\) on this interval are at least
two by unimodality and the endpoint values at \(1,a-2\).
Thus \(q_{X,K-2}-1\ge1=C_0\) here as well.
\end{proof}

\begin{lemma}[Reverse half-range]\label{lem:half-tail}
For the reversed product of two increasing arms,
\begin{equation}\label{eq:half-tail}
 T_p(n)\le0\qquad\bigl(n>\lastmode(p),\ 2n\le a+b\bigr).
\end{equation}
\end{lemma}
\begin{proof}
For a uniform window \(w=[t]p\),
\(T_w(n)=E_p(n+1)-E_p(k)\), where \(k=n+1-t\).
If \(n>\lastmode(w)\) and \(k\le\lastmode(p)\), then
\[
 \Delta w_{n+1}-\Delta w_n
 =\Delta p_{n+1}-\Delta p_k\le0 .
\]
Indeed, write \(m=\firstmode(p)\). If \(m>0\), then
\(\Delta w_m=p_m-p_{m-t}>0\), so \(\firstmode(w)\ge m\). The same
mode bound is immediate for \(m=0\). Thus \(\Delta p_{n+1}\le0\),
whereas \(\Delta p_k\ge0\) for \(k\le\lastmode(p)\), also for \(k<0\)
by zero extension. Since \(\Delta w_{n+1}\le0\), we have \(T_w(n)\le0\).

For the base with an interval factor, write \(p=[a+1]q_Y\).
If \(a\ge b\) the maximum interval is \([b,a]\), so the range in \eqref{eq:half-tail} is empty. If both factors are intervals, the product is symmetric
and unimodal, so its last mode is at or beyond its midpoint. Otherwise \(a<b\) and
\(q_Y\) is a noninterval arm. Below rank \(b\),
\(2q_{Y,j}=(R_Y)_j+q_{Y,j-1}\).
Through the last mode of \(R_Y\), this averaging recurrence
gives strict increase of \(q_Y\) because its preceding
value is strictly less than the current positive source.
Thus \(E_{q_Y}\) is nonincreasing after the ordinary mode
and below \(b\). Apply the window observation above and,
when \(k>\lastmode(q_Y)\), compare the two source indices
in this nonincreasing range. The endpoint \(n+1=b\) can occur only for \(b-a\in\{1,2\}\).
By \eqref{eq:sources},
\[
 T_p(b-1)=1-(R_Y)_{b-a-1}\le0,
\]
since \((R_Y)_0=2\) and \((R_Y)_1\ge1\) for a
noninterval arm.

Induct on the total number of parts greater than one after normalization.
Orient the arms so that \(a\ge b\), and write \(X=Zt\). With
\(v=q_Zq_Y\), we have \(p=[t]v+z^a q_Y\). Its first mode lies
below \(a\) by Lemma~\ref{lem:ordinary}. Up to that mode it agrees with
the window. The window is bounded above by \(p\) everywhere,
so \(\lastmode([t]v)\le\lastmode(p)\).
If the term \(z^a q_Y\) can affect \(T_p(n)\), then \(a\le n+1\).
Since \(2n\le a+b\le n+1+b\), we have \(n\le b+1\). The forward index
\(j=a+b-n\) therefore satisfies \(j+1\ge a\ge b\). Both forward factors
are nonincreasing at that edge, since each has top coefficient one
and positive coefficients throughout its support. Also
\(2n\ge2a-2\ge a+b-2\), giving \(2j\le a+b+2\).
Lemma~\ref{lem:both-nonrise} now yields \(T_p(n)\le0\).

On every other edge, \(n+1<a\), so \(p\) agrees with \([t]v\)
at all indices used in \(T_p(n)\), and \(k=n+1-t<a-t\).
The window calculation at the start of the proof gives
\(T_{[t]v}(n)\le0\) whenever \(k\le\lastmode(v)\).
If the shortened arm \(Z\) is an interval arm, the last mode of
\(v\) is at least \(a-t=\wt(Z)\). Indeed its differences
through \(a-t\) are the positive coefficients of \(q_Y\) as long as
they are in the support and are zero thereafter. Thus \(k<a-t\) lies at or before that mode.

It suffices to consider two remaining noninterval arms with
\(k>\lastmode(v)\). Call the heavier arm \(L\) and the lighter one \(S\),
breaking ties arbitrarily. Write \(d=\wt(L)\), \(s=\wt(S)\),
\(K=\kappa(L)\), \(R=d+1-K\), and \(C=R_L q_S\). By \eqref{eq:sources},
\[
 E=(2-z)v=C+z^d(1-z)q_S .
\]
Below rank \(d\) the correction vanishes, so \(E=C\).
Each interval-product component
of \(C\) has endpoint sum at most \(d+s-1\). Consequently
\(C\) is nonincreasing from \(\lfloor(d+s)/2\rfloor\)
onward. For \(i>\lastmode(v)\) with \(2i\le d+s\), the induction hypothesis
gives \(E_{i+1}\le E_i\). Together with the bound for \(C\), this
makes \(E\) nonincreasing throughout \(\lastmode(v)<i<d\).
Since \(k<a-t\le d\), for \(n+1<d\) we obtain
\(T_{[t]v}(n)=E_{n+1}-E_k\le0\).

At \(n+1=d\), the correction adds one, giving \(E_d=C_d+1\).
If \(K\le s+2\) and \(d>s\), the coefficient in the product of
\([0,R]\) with \([0,s]\) falls by one between ranks \(d-1\) and \(d\),
since the leaving point \(K-2\) lies in \([0,s]\)
and the new entering value is zero.
When \(d=s\), the second height interval of \(q_S\) is
\([1,s-2]\). It contains \(K-2\) because \(K\ge3\) and
\eqref{eq:kappa} bounds \(K-2\) by the number of parts greater than
one in \(L\), at most \(\lfloor(s-1)/2\rfloor\).
All other products have nonnegative drops. Thus
\(E_d\le E_{d-1}\le E_k\).
If \(K>s+2\), then \(C_d=0\) and \(E_d=1\).
The hypothesis \(2n\le a+b\) gives \(k=d-t\le s+2\),
whereas Lemma~\ref{lem:head-range} gives \(R\ge2\)
for \(K\ge s+3\ge6\). The support of \(C\) therefore
includes \(k\), so \(E_k=C_k\ge1\).

Above \(d\), write \(n=d+r\) with \(r\ge0\). The assumption
\(2n\le a+b\) gives \(k\le k_*=s-r+1\). The correction at \(n+1\)
is now a first difference of \(q_S\), so it suffices to prove
\begin{equation}\label{eq:half-tail-head}
 C_{d+r+1}+q_{S,r+1}-q_{S,r}\le E_k .
\end{equation}
Every interval product of \(C\) has a nonnegative drop
from \(d-1\) to \(d+r+1\) because the two indices sum
to \(2d+r\).
If \(r+3\ge K\), pair each ordinary layer \([A,B]\)
starting at \(A=r+1\) with the lowest source layer \([0,R]\).
At \(d-1\) it counts the whole layer, since \(B\le s-r-1\)
and \(d-1-R=K-2\le A\). At \(d+r+1\) it omits \(A\).
There are at least \(q_{S,r+1}-q_{S,r}\) such layers, so
their total drop proves \eqref{eq:half-tail-head} using
\(E_{d-1}\le E_k\).

For \(0\le r\le K-4\), write \(J_S\) for the main interval product of
the lighter arm. Its degree is \(s-h_S\), where \(h_S\) counts the
parts of \(S\) greater than one. If
\(\deg J_S\ge K+r-1\), the drop of the lowest source interval
convolved with \(q_S\), from \(d-1\) to \(d+r+1\), is
\[
 \sum_{u=K-2}^{K+r-1}q_{S,u}-\sum_{u=d}^{d+r+1}q_{S,u}.
\]
The second sum is at most one since \(d\ge s\). We have
\(r+1\le K-2\) and \((r+1)+(K-2)\le\deg J_S\), so symmetry
and unimodality give \((J_S)_{K-2}\ge(J_S)_{r+1}\).
Lemma~\ref{lem:head-range} then gives
\(q_{S,K-2}\ge(\Delta q_{S,r+1})_+\). The term
\(q_{S,K-1}\ge1\) covers the possible subtraction. Thus the drop
is at least \(\Delta q_{S,r+1}\), and \(E_{d-1}\le E_k\) gives
\eqref{eq:half-tail-head}.

If \(k_*<d\) and \(K\le d-s+2r+1\), compare instead
with \(C_{k_*}\). The two indices sum to \(d+s+2\),
so all components have nonnegative drop.
A layer of \(q_S\) beginning at \(r+1\) lies wholly
in the convolution with \([0,R]\) at \(k_*\), since
its upper endpoint is at most \(s-r-1\), and
\(k_*-R\le r+1\). The window at \(d+r+1\) omits its
starting term. This proves \eqref{eq:half-tail-head}
because \(C_{k_*}=E_{k_*}\le E_k\).
These two comparisons exhaust \(k_*<d\). Indeed, if both failed, then
\(d\le s+K-2r-2\) and \(\deg J_S\le K+r-2\).
Normalization gives \(2\deg J_S\ge s+1\), hence \(d\le3K-7\),
contrary to \(d\ge K(K-1)/2>3K-7\).

If \(k_*=s-r+1\ge d\), then \(r+d-s\le1\). Thus \(r\le1\), with
\(d=s\) when \(r=1\). For \(r=0\) the same weight bound
gives \(d\ge2K-3\), so \(2\deg J_S\ge s+1\ge d\) implies
\(\deg J_S\ge K-1\). For \(r=1\), the bound \(r\le K-4\) gives
\(K\ge5\). Then \(d=s\ge K(K-1)/2\ge2K-1\), so \(\deg J_S\ge K\).
In both cases the drop from the lowest source interval proves
\eqref{eq:half-tail-head}.
\end{proof}

Together the two lemmas cover the case of two nonincreasing forward
factors. Lemma~\ref{lem:both-nonrise} applies when \(2j\le a+b+2\). Otherwise
\(2n<a+b-2\), and \(n>\lastmode(p)\) in Condition~\ref{cond:local},
so Lemma~\ref{lem:half-tail} again gives \(T_p(n)\le0\).
\subsection{When one factor is still increasing}

Normalize both arms and use \(s_Y\) from Lemma~\ref{lem:head-range}.
Suppose \(B\) rises from \(j\) to \(j+1\) while \(A\) does not.
Then the increments of \(\mathcal E_Y\) contribute with one sign, and
the term they must offset is a second difference of \(q_X\).

\begin{proposition}[The contribution from the nonincreasing arm]\label{prop:mixed-head}
Under the hypotheses of Condition~\ref{cond:local}, orient the
arms so that \(A_{j+1}\le A_j\) and \(B_{j+1}>B_j\), where \(j=a+b-n\).
Put \(r=a-j\). The comparison \eqref{eq:local} holds whenever
\(T\le0\), \(r<0\), \(r\ge k\), or \(j-2\le s_Y\).
It also holds if the normalized \(X\) has only one part greater than one.
\end{proposition}
\begin{proof}
The right side of \eqref{eq:local} is nonnegative because
\(k<\ell\), so assume \(T>0\). We first locate \(k\) relative to \(j\). The reversed
interval products have centers at most \((a+b)/2\),
so \(p_i\ge p_{i'}\) for \(i\le i'\), \(i+i'\ge a+b\).
Since \(p_{n+1}>p_{k+1}\) with \(n>k\), we have
\(n+k+2<a+b\), or equivalently \(k+1\le j-2\).

Let \(v_i=(\Delta q_{X,i})_+\) and
\(C_0=\Delta q_{X,r+1}-\Delta q_{X,r}\). Reversing the
identity for \(e_i\) in Lemma~\ref{lem:tilt} gives
\begin{equation}\label{eq:mixed}
 T=C_0-\sum_{i=2}^{j+1}(\Delta\mathcal E_Y)_i A_{j+1-i},
 \qquad E_n\ge \Delta q_{X,r}.
\end{equation}
The source has zero coefficients at ranks zero and one, and terms
above \(j+1\) have \(A_{j+1-i}=0\). Through rank \(j+1\) the source
increments are nonnegative by the argument of Lemma~\ref{lem:tilt},
since \(B\) is still rising. Hence
\(C_0\ge T>0\). Our orientation gives \(\Delta q_{X,r}\ge0\).
For \(r\ge0\), the positive support of \(v\) is initial, so
\(v_{r+1}>v_r>0\). This rules out an arm \(X\) with just one part
greater than one at these ranks. Its main product
\([\alpha+1][t]\) is trapezoidal, and the terminal unit adds no
positive difference.

Suppose \(r\ge k\). Set \(\lambda=C_0/v_r\). By
Lemma~\ref{lem:rising-lc},
\(v_{i+1}-v_i\ge\lambda v_i\) for \(0\le i\le k\).
Telescoping gives \(\lambda q_{X,i}\le v_{i+1}-v_0\)
for \(i\le k\). Add the inequalities to get
\[
 2v_{i+1}-v_i\ge\lambda(q_{X,i}+v_i)+1 .
\]
The expressions \(2v_{i+1}-v_i\) and \(q_{X,i}+v_i\) are the
difference at \(i+1\) and the coefficient at \(i\) of \((2-z)q_X\).
Multiply the inequality by \(q_{Y,k-i}\) and sum for \(0\le i\le k\).
The constants contribute \(\sum_{i=0}^kq_{Y,k-i}\ge q_{Y,0}=1\),
and the remaining convolution term is \(2q_{Y,k+1}\ge0\). Thus
\(E_{k+1}-E_k\ge\lambda E_k+1\).
Since \(E_k\ge E_n\ge v_r\), this is at least \(C_0+1\).

For \(r<k\) and \(j-2\le s_Y\), use instead
\[
 E_{k+1}-E_k=\coef{k+1}{(2-z)q_X(1-z)q_Y}.
\]
One has \(r\ge-1\) since \(C_0=0\) for \(r\le-2\).
Consequently \(k+1\le j-2\le a-1\).
By \eqref{eq:sources}, the coefficients of \((2-z)q_X\) through
degree \(k+1\) agree with those of \(R_X\) and are nonnegative.
The differences of \(q_Y\) through that degree are all at least one. If \(r\ge0\), the term
at index \(r+1\) is at least
\(2q_{X,r+1}-q_{X,r}\ge C_0\).
If \(r=-1\), the constant coefficient of \((2-z)q_X\) is two and
\(C_0=1\). Thus in both cases \(E_{k+1}-E_k\ge C_0\ge T\).

Finally, for \(r=-1\) one has \(j=a+1\).
If \(K=\kappa(Y)\le j+1\), the first source term alone
gives \(A_{j+1-K}\ge1=C_0\), because \(K\ge3\)
for a noninterval arm. Otherwise \(j\le K-2\), so
\(j-2\le K-4\le S_K\le s_Y\). The previous bound then applies to \(r=-1\) as well.
\end{proof}

Let \(K=\kappa(Y)\). If the largest part of \(X\) is at least \(K-1\)
and \(C_0>0\), Lemma~\ref{lem:capacity} gives
\(C_0\le q_{X,r+K-1}\). This coefficient is positive, so \(K\le j+1\)
and the first nonzero source increment occurs in \eqref{eq:mixed}.
It contributes at least \(q_{X,r+K-1}\), giving \(T\le0\).
It remains to prove the comparison for \(0\le r<k\) and \(j-2>s_Y\)
when the largest part of \(X\) is at most \(K-2\).
In this range a bound on \(\Delta B_{j+1}\) suffices for
\(K\ge13\). For \(K\le12\) the finite reduction is
Proposition~\ref{prop:small-k}.
\begin{lemma}[Rank growth in the reversed arm]\label{lem:mark-growth}
For a normalized increasing arm \(X=1^\alpha(t_1,\ldots,t_h)\)
with \(\alpha\ge1\), write \(Q_i=q_{X,i}\). Then
\begin{equation}\label{eq:mark-growth}
 (r+1)Q_{r+1}\le(r+h+1)Q_r \quad(r\ge0).
\end{equation}
In particular \(Q_1=h+1\).
\end{lemma}
\begin{proof}
Start with a chain of \(\alpha\) elements. For each successive part
\(t\), take the disjoint union of the current poset with a chain of \(t-1\) elements and
adjoin a new maximum above both. An ideal containing that maximum is the entire new poset. Ideals
omitting it can choose an initial segment of the new chain independently
of the ideal in the old poset. Their generating function thus obeys the
recurrence for \(q_X\) in \eqref{eq:append}. The cover graph is a tree rooted at
its maximum, with \(h+1\) minimal elements.

Take an ideal of size \(r+1\) and a marked element. Follow
the unique upward path from the marked element as far as it
stays in the ideal and delete its last element \(v\).
This leaves an ideal \(I\) of size \(r\). For fixed \(I\),
the possible \(v\)'s are minimal elements of the complement.
They are incomparable, so their principal lower sets are
disjoint. Each contains a minimal element of the full
poset, so there are at most \(h+1\) such \(v\)'s.
For a fixed \(v\), the mark must belong to its principal lower set.
Their disjoint union contains at most the \(r\) elements of \(I\) and
\(h+1\) choices of \(v\). Counting the marked ideals in these two
ways gives \eqref{eq:mark-growth}.
\end{proof}

\begin{lemma}[Two bounds on late forward slopes]\label{lem:dyadic-slope}
For a normalized increasing arm \(Y\), put \(B=F_Y\) and retain \(K=\kappa(Y)\)
and \(M_Y,B_Y,s_Y\) from Lemma~\ref{lem:head-range}.
At any integer \(m>s_Y\),
\begin{align}
 20\,\Delta B_m&\le 2^m &&(K\ge13),\label{eq:constant-density}\\
 m\,\Delta B_m&\le 2^m &&(K\ge16).\label{eq:rank-density}
\end{align}
\end{lemma}
See Appendix~\ref{app:slope} for the proof.

\begin{proposition}[Large first missing weights]\label{prop:large-k}
Under the hypotheses of Condition~\ref{cond:local}, suppose
\(A_{j+1}\le A_j\) and \(B_{j+1}>B_j\) at \(j=a+b-n\).
If \(K=\kappa(Y)\ge13\), then \eqref{eq:local} holds.
\end{proposition}
\begin{proof}
By Proposition~\ref{prop:mixed-head} it remains to treat
\(T>0\), \(0\le r<k\), and \(j-2>s_Y\), where \(r=a-j\) and
the normalized arm \(X\) has \(h\ge2\) parts greater than one.
Set \(N=j+1\),
\(Q=q_{X,r+1}\), and \(R=2q_{X,r}-q_{X,r-1}\ge q_{X,r}\).
Then \(C_0=Q-R\) in \eqref{eq:mixed}.
For \(2\le i\le N\) the increments \(w_i=(\Delta\mathcal E_Y)_i\)
are nonnegative. Since \(A_{N-i}\ge2^{2-i}A_{N-2}\)
for \(2\le i\le N\) by Lemma~\ref{lem:count-bounds},
and \(A_{N-2}=Q\), the telescoping identity
\[
 4\sum_{i=2}^N2^{-i}w_i=1-4\,2^{-N}\Delta B_N
\]
gives the lower bound \(Q(1-4\,2^{-N}\Delta B_N)\)
on the source sum. Here the coefficients \(2B_{i-1}-B_i\) for \(i\ge2\)
telescope, with \((\mathcal E_Y)_1=0\) and \(B_0=B_1=1\).

For \(K\ge16\), we claim \(NR\ge4Q\).
The inequalities \(j\ge r+4\) and \(j+r=a\ge2h+1\)
follow from \(k+1\le j-2\) and normalization.
For \(r\ge3\) set \(x=r+1\). If \(h\le x\), then
\(Nx\ge(x+4)x\ge4(h+x)\).
If \(h\ge x+1\), then
\(Nx-4(h+x)\ge(2x-4)h-x(x+1)\ge (x-4)(x+1)\ge0\).
Thus the claim follows from \eqref{eq:mark-growth}.
When \(r=2\), \(3N\ge6h\ge4(h+3)\) if \(h\ge6\).
If instead \(h\le5\) at this rank, \(j-2>s_Y\ge S_K\) gives \(N\ge11\), again enough
for \(3N\ge4(h+3)\). At \(r=1\) we have \(R=2h+1\),
\(2Q\le(h+1)(h+2)\), and \(N\ge2h+1\). Thus
\(NR-4Q\ge2h^2-2h-3>0\). At \(r=0\) the claim follows from
\(R=2\), \(Q=h+1\), and \(N\ge2h+2\).
Since \(N>s_Y\), \eqref{eq:rank-density} bounds the source sum
below by \(Q-4Q/N\ge Q-R=C_0\), contrary to \(T>0\).

For \(13\le K\le15\), the main interval product \(J_X\) has
degree \(a-h\). If \(\deg J_X\ge2r+K\), its symmetry and
unimodality give
\(q_{X,r+K-1}\ge(J_X)_{r+1}\ge C_0\), by
Lemma~\ref{lem:head-range}. The first term of the forward
source sum already gives this coefficient. Therefore
\(\deg J_X\le2r+K-1\).
If \(h>4(r+1)\) and \(r\ge5\), then
\(\deg J_X\ge h+1\ge4r+6\ge2r+K\), a contradiction.
If \(h>4(r+1)\) and \(r\le4\), then
\(j-2=\deg J_X+h-r-2\le3r+2K-5\le S_K\)
for \(K=13,14,15\), also a contradiction.
Hence \(h\le4(r+1)\) and
\eqref{eq:mark-growth} gives \(Q\le5q_{X,r}\le5R\).
By \eqref{eq:constant-density} the source sum is at least
\(4Q/5\ge Q-R=C_0\), which contradicts \(T>0\) in this range also.
\end{proof}

\subsection{Proof of the weighted difference bound}

\begin{proposition}[Weighted positive differences]\label{prop:count-induction}
For every pair of increasing arms the sequence \(U\) is unimodal.
Consequently,
\[
 \sum_{i\ge1}u_{k+i}\le u_k\qquad(k\ge\firstmode(U)).
\]
\end{proposition}
\begin{proof}
Normalize the arms and induct on the total number of parts greater
than one. Lemma~\ref{lem:weighted-bases} covers an interval factor and
a heavier arm with a single part greater than one. Otherwise write the
heavier arm as \(X=Zt\). Both \(Z\) and \(Y\) are noninterval, and
the weighted positive differences of \(q_Zq_Y\) are unimodal by induction.

At a rank where both forward factors of the shorter product rise,
Lemma~\ref{lem:tilt} gives the inequality in Condition~\ref{cond:local}.
If both factors are
nonincreasing, Lemmas~\ref{lem:both-nonrise} and \ref{lem:half-tail}
give the inequality on the two halves of the rank range. When just one
factor rises, the comparison follows from Proposition~\ref{prop:large-k}
for \(K\ge13\) and Proposition~\ref{prop:small-k} for \(K\le12\).
Lemmas~\ref{lem:ordinary} and \ref{lem:tilt} supply the other hypotheses
of Lemma~\ref{lem:window}. Hence the weighted positive differences of
\([t]q_Zq_Y\) are unimodal.

By Lemma~\ref{lem:ordinary}, \(p=[t]q_Zq_Y+z^a q_Y\) reaches its
first maximum before the added term begins. Its weighted positive
differences are therefore the initial portion of those of \([t]q_Zq_Y\),
followed by zeros. Finally, for \(k\ge\firstmode(U)\) we have
\(u_{k+i}\le2^{-i}u_k\). The bound on their sum follows from
\(\sum_{i\ge1}2^{-i}=1\).
\end{proof}
\section{Unimodality at a strict peak}\label{sec:peaks}

Let \(X,Y\) be increasing arms and let \(M\) be a positive integer.
We call \(M\) a strict peak if \(M>\max\{\max(X),\max(Y)\}\).
For such a peak, set
\[
 H_M(X,Y)=F_X+F_Y+(S_M-1)F_XF_Y,\qquad
 Q_M(X,Y)=F_X\{1+(S_M-1)F_Y\}.
\]
By \eqref{eq:core}, \(H_M(X,Y)\) counts the subwords of \(XMY\).
We call \(Q_M(X,Y)=H_M(X,Y)-F_Y\) the relative polynomial.
The deletion identities \eqref{eq:relative-append} imply that in a
relative counterexample with the fewest arm parts at fixed \(M\),
the arm weights \(a=\wt(X)\) and \(b=\wt(Y)\) satisfy
\(\lvert a-b\rvert\le M-1\). In this range \(H_M\) and both relative
polynomials are unimodal under the same conditions
(Lemma~\ref{lem:balanced}). We can therefore work with \(H_M\),
whose middle-rank differences involve the product bound of
Proposition~\ref{prop:count-induction}.

\begin{theorem}[Peak reduction]\label{thm:peak}
The polynomials \(H_M(X,Y)\), \(Q_M(X,Y)\), and \(Q_M(Y,X)\)
have unimodal coefficients at every strict peak.
\end{theorem}
\begin{proof}
Suppose that a relative polynomial \(Q_M(X,Y)\) fails to be unimodal.
Choose the smallest such \(M\), then the fewest parts in the two arms
together. By Lemmas~\ref{lem:balanced} and \ref{lem:outer} the arms are
noninterval and we can orient them with \(a=\wt(X)\le b=\wt(Y)\)
and \(b-a\le M-4\). It suffices to prove unimodality of the core
\(H_M\) for this pair. Write \(P=F_XF_Y\) and \(p=q_Xq_Y\).
Put \(N=a+b+M\), \(R_0=\lfloor N/2\rfloor\), and
\[
 W=([M]-z^M)p,\qquad
 \delta_M(Y)=([M]-z^M)q_Y+z^{b+M}.
\]
The reversal of \(Q_M(X,Y)\) about degree \(N\) is
\(q_X\delta_M(Y)\). The reverse core equals \(W\) through
rank \(R_0\) because its two corrections begin at \(a+M\)
and \(b+M\), both larger than \(R_0\).

Write \(Y=Zt\), \(c=b-t\).
By minimality, the reversal \(q_X\delta_M(Z)\) of the relative
polynomial with the shorter arm \(Z\) is unimodal. Its convolution with \([t]\) agrees with \(W\)
through \(b-1\), since
\begin{equation}\label{eq:proper-relative}
 W-[t]q_X\delta_M(Z)
 =z^b q_X([M]-z^M)-z^{c+M}[t]q_X ,
\end{equation}
and \(c+M>b\). Up to \(b-1\), therefore, the differences of \(W\)
can change sign only from positive to negative. For \(n<M\) we also have
\(\Delta W_n=p_n\ge0\), since \((1-z)W=(1-2z^M+z^{M+1})p\).
It remains to determine the signs starting at \(\max(b,M)\).

Put \(\sigma=\firstmode(p)<b\) and
\(\ell=\firstmode(U)\), with \(u,U\) as in
Section~\ref{sec:product}. For
\(\max(b,M)\le n=M+k\le\min(R_0,M+\sigma)\), we have
\begin{equation}\label{eq:joining}
 \Delta W_{M+k}=p_{M+k}-p_k-u_k .
\end{equation}
Indeed \(k\le\sigma\) gives \(\Delta p_k=u_k\). While \(k<\ell\),
\(p_k+u_k\) increases weakly since \(2u_{k+1}-u_k\ge0\).
Over the same ranks \(p_{M+k}\) decreases weakly, as \(M+k\ge b>\sigma\).
If \(b>M\), this comparison extends one rank back to \(b-1\), since
\(b-1\ge\sigma\) and \(k=b-1-M\ge0\) when the range is nonempty.
Thus a negative difference already present in the prefix cannot be
followed by a positive difference in the portion with \(k<\ell\).
If \(b\le M\), the preceding differences are nonnegative.
Once \(\ell\le k\le\sigma\), Proposition~\ref{prop:count-induction}
gives
\[
 p_{M+k}-p_k\le\sum_{i=1}^{M}u_{k+i}\le u_k.
\]
Substitution in \eqref{eq:joining} gives \(\Delta W_{M+k}\le0\)
through \(k=\sigma\) whenever that rank lies in the stated range.

For any remaining rank \(M+\sigma<n\le R_0\), put
\(s=a+b-\sigma=\lastmode(P)\) and \(j_H=N-n<s\).
If this range is nonempty, then \(s>(a+b+M)/2>b\). Lower the peak
to the largest part \(t_0\) in either arm. At peak \(t_0<M\), remove
the end copies of \(t_0\) from both arms. The remaining peak is strict.
Its relative polynomials are unimodal by minimality of \(M\), and
Lemma~\ref{lem:balanced}, with the heavier arm first, gives unimodality
of its core. With both arms empty the core is simply an interval.
Restoring a copy of \(t_0\) replaces its arm polynomial \(F_Z\) by
\(1+S_{t_0}F_Z\). Substitution in \eqref{eq:core} gives
\(H_{t_0}(Zt_0,Y)=1+S_{t_0}H_{t_0}(Z,Y)\), with the symmetric
identity for the other arm. The window is unimodal and begins at
rank one with coefficient one, so adding the constant term preserves
unimodality. Restoring all copies gives unimodality of \(H_{t_0}(X,Y)\).

As \(s>b\ge a\), the arm differences vanish at \(s+1\).
The product has its last forward mode at \(s\), so
\[
 \Delta (H_{t_0})_{s+1}
 =2P_s-P_{s+1}-P_{s-t_0}>0 .
\]
Hence \(H_{t_0}\) is nondecreasing through \(s+1\). From
\(H_M=H_{t_0}+z^{t_0+1}[M-t_0]P\), for the corresponding
\(j_H<s\) we obtain
\[
 \Delta (H_M)_{j_H+1}=\Delta (H_{t_0})_{j_H+1}+P_{j_H-t_0}-P_{j_H-M}\ge0 ,
\]
since \(P\) is nondecreasing through \(s\). Reversing the indices gives
nonpositive differences through \(R_0\) on this last interval. Beyond
\(R_0\) they are nonpositive by Lemma~\ref{lem:count-bounds} applied
to \(XMY\). Hence \(H_M\) is unimodal, contrary to the equivalence
with \(Q_M\) for our pair of arms. All relative polynomials are
therefore unimodal. For an arbitrary strict peak, taking the heavier
arm first in Lemma~\ref{lem:balanced} gives unimodality of \(H_M\)
as well.
\end{proof}

\begin{proof}[Proof of Theorem~\ref{thm:main}]
First suppose that \(w\) has exactly one descent, and let \(M\) be
the part immediately before it. Let \(r+1\) be the length of the
constant block of \(M\)'s ending at the descent. Then \(w=XM^{r+1}YE\),
where \(X,Y\) are increasing with all parts less than \(M\), and the
suffix \(E\) is increasing with all parts at least \(M\).
Theorem~\ref{thm:peak} gives unimodality for \(XMY\).
Each additional central copy of \(M\) sends its rank polynomial
\(F\) to \(1+S_MF\), as in the proof above. Each successive part
\(t\) of \(E\) is a new or repeated maximum and sends \(F\) to
\(1+S_tF\). These operations preserve unimodality by
Lemma~\ref{lem:layers}, since \(S_tF\) starts with coefficient one
at rank one. For a word without descents the same append argument
starts at the empty word.

For a nonempty bitonic word, retain one central maximum and remove
its other copies. Reversing the decreasing right arm preserves its
rank polynomial. Thus \eqref{eq:core} gives the same rank polynomial
as for a strict peak with two increasing arms, allowing either arm to
be empty. Theorem~\ref{thm:peak} applies to this peak, and restoring
the repeated maxima preserves unimodality as before.
\end{proof}

\section*{Statements and Declarations}
\addcontentsline{toc}{section}{Statements and Declarations}
\noindent\textbf{Funding.} No funding was received for this work.

\noindent\textbf{Competing interests.} The authors declare that they have no competing interests.

\noindent\textbf{Author contributions.} The authors carried out the research jointly and contributed equally.

\noindent\textbf{Data and code availability.} Online Resource~1 contains the two exact
Python programs, the critical inequalities, and instructions for reproducing
Table~\ref{tab:finite}. The earlier arithmetic certificates and Lean proofs of
selected auxiliary statements are archived at
\href{https://github.com/bos024/one-descent-rank-unimodality/tree/b0840c7155fe2db0716061c8970e1041311349a8}{GitHub revision \texttt{b0840c7155fe}}.

\noindent\textbf{Online Resource~1.} Code, exact verification data, and instructions for reproducing Table~\ref{tab:finite}, together with the Lean source and log for four auxiliary integer inequalities.

\noindent\textbf{Use of generative AI.} In the original research, OpenAI Codex
was used solely to assist with the Lean formalization of selected results.

\clearpage
\begingroup\singlespacing\small
\renewcommand{\bibsection}{\section*{\refname}%
  \addcontentsline{toc}{section}{\refname}}

\endgroup

\clearpage
\appendix
\setstretch{1.25}
\section{Relative polynomials and the outer cases}\label{app:peak}

In a minimal counterexample at a fixed peak \(M\), the two arm weights
must be close. We obtain the bound \(\lvert a-b\rvert\le M-4\) from
the two ways of deleting a part in a relative polynomial.

For the increasing arms \(X,Y\), let \(a=\wt(X)\), \(b=\wt(Y)\),
\(B=F_Y\), and
\[
 K_M(Y)=1+(S_M-1)B,\qquad R=b+M.
\]
Then \(Q_M(X,Y)=F_X K_M(Y)=H_M(X,Y)-F_Y\), and
\(\delta_M(Y)\) in Section~\ref{sec:peaks} is the reversal of
\(K_M(Y)\) about degree \(R\).
At a strict peak, applying \eqref{eq:append} to each arm in turn gives
\begin{equation}\label{eq:relative-append}
 Q_M(Zu,Y)=K_M(Y)+S_u Q_M(Z,Y),\qquad
 Q_M(X,Yv)=S_v Q_M(X,Y)+(S_M-S_v)F_X .
\end{equation}

\begin{lemma}[Relative kernel]\label{lem:kernel}
The polynomial \(K_M(Y)\) is nonnegative and unimodal of
degree \(R\), and all its height intervals have endpoint
sum at least \(R+1\). For nonempty \(Y\) it has zero
coefficients at \(0,1\) and coefficient one at \(R\).
Consequently \(Q_M(X,Y)\) is nondecreasing through rank
\(\lfloor(a+b+M+2)/2\rfloor\).
\end{lemma}
\begin{proof}
The empty arm gives \(S_M\), with a single height interval \([1,M]\).
If \(Y=Zt\), then \(K_M(Y)=S_tK_M(Z)+S_M-S_t\). Since \(t<M\),
the recurrence preserves nonnegativity. Also
\[
 (1-z)K_M(Z)=\mathcal E_Z-z^{M+1}F_Z
\]
by \eqref{eq:sources}, so \(K_M(Y)\) increases through rank
\(M\). Write \(P=K_M(Z)\) and \(V=S_tP\). At rank \(M+1\) the difference is \(P_M-P_{M-t}-1\).
If \(P_M>P_{M-t}\), this difference is nonnegative by integrality.
All subsequent coefficients belong to the unimodal window. Otherwise \(\sum_{i=M-t+1}^M(\mathcal E_Z)_i=0\).
By Lemma~\ref{lem:source}, its support ends at or before
\(M-t\), since its initial index is at most \(t+1\le M\).
The sequence \(P\) is constant as long as these differences
vanish and nonincreasing thereafter. Thus all subsequent
window differences are nonpositive. This proves unimodality.

For the interval claim, write \(c=\wt(Z)\), \(d=c+M\).
Inductively \(P\) has endpoint sums at least \(d+1\).
Every interval product in \(V=S_tP\) then has endpoint sum
at least \(d+t+2\).
In comparing coefficients of \(K_M(Y)\) at \(i\) and its
reflection \(d+t+1-i\), the added interval \([t+1,M]\)
can cause a loss only if
\[
 t+1\le i\le M,\qquad i\le c+t,\qquad 2i<d+t+1.
\]
Set \(f_j=P_{d+1-j}-P_j\). The difference from the window is
\(\sum_{j=i-t+1}^i f_j+P_i-P_{i-t}\).
Any negative \(f_j\) in this window has its nonnegative
opposite \(f_{d+1-j}\) in the same window, because the sum
of its endpoints is at most \(d+1\). The last difference
is a sum of source coefficients, since \(i\le M\).
The window meets the source support \([k,c+1]\) by
\(k\le t+1\) and \(i\le c+t\), so it contributes at least one. It follows that every left-half coefficient
of \(K_M(Y)\) is at most its reflection about \((R+1)/2\).
For each height interval, apply this comparison to its left endpoint
if that endpoint is at most \((R+1)/2\). Its right endpoint must
reach the reflection. A height interval starting past the center
already has endpoint sum greater than \(R+1\). This proves the
interval bound. For nonempty \(Y\), the two initial zeros and final
coefficient one follow from \(1+(S_M-1)B\). Taking the products with
the height intervals of \(F_X\) now gives the claimed increasing range
for \(Q_M\).
\end{proof}

\begin{lemma}[Strict reflected difference]\label{lem:tail-budget}
For \(B=F_Y\) and \(1\le v\le\min(b-1,M-1)\),
\begin{equation}\label{eq:reflected-lower}
 (K_M(Y))_{R-v}-(K_M(Y))_{v+1}
 =\sum_{j=0}^v(B_{b-j}-B_j)+B_{v+1}\ \ge B_{v+1}\ge1.
\end{equation}
\end{lemma}
\begin{proof}
At rank \(R-v\) the subtracted \(B\) term is absent and the
window sums the last \(v+1\) coefficients of \(B\).
At \(v+1\) it sums the first \(v+1\) and subtracts
\(B_{v+1}\). The reflected differences in the sum are
nonnegative before the midpoint by Lemma~\ref{lem:layers}. Pairing
opposite indices past the midpoint leaves only nonnegative terms. The support of \(B\) contains \(v+1\).
\end{proof}

\begin{lemma}[Balanced deletion]\label{lem:balanced}
For a fixed \(M\), any nonunimodal relative with the fewest
total arm parts has nonempty noninterval arms with
\(|a-b|\le M-1\). More generally, if \(b\le a+M-1\),
\(Q_M(X,Y)\) is unimodal if and only if \(H_M(X,Y)\) is
unimodal. In the balanced range the arms can be exchanged.
\end{lemma}
\begin{proof}
Empty arms give \(K_M(Y)\) or \(S_M F_X\), both unimodal.
If \(a\ge b+M\), the first recurrence in
\eqref{eq:relative-append} has its correction vanish after
degree \(b+M\). Lemma~\ref{lem:kernel} gives increase through
rank \(b+M+1\), from which point the coefficients agree with
a window of the relative polynomial with the shorter arm, unimodal
by minimality.
If \(b\ge a+M\), the correction in the second recurrence has degree at
most \(a+M\), and \(Q_M\) increases through \(a+M+1\). A counterexample
of minimal length must therefore have \(|a-b|\le M-1\).

If \(b\le a+M-1\), both \(H_M\) and \(Q_M(X,Y)\) increase
through \(b+1\), by Lemmas~\ref{lem:count-bounds} and
\ref{lem:kernel}. Their later differences agree since
\(H_M-Q_M(X,Y)=B\). This gives the equivalence, in both orientations when \(|a-b|\le M-1\).
If \(F_X\) is an interval then \(F_XK_M(Y)\) is unimodal by
Lemma~\ref{lem:layers}. Exchanging the arms deals with an interval \(F_Y\).
\end{proof}

\begin{lemma}[The outer three cases]\label{lem:outer}
For \(M\ge3\), a least relative counterexample as in
Lemma~\ref{lem:balanced} satisfies \(|a-b|\le M-4\).
For \(M\le2\) there are no counterexamples.
\end{lemma}
\begin{proof}
For \(M\le2\) all possible arms have interval polynomials.
Orient a least counterexample with \(a\ge b\).
Write \(X=Zu\), \(R=b+M\), \(V=F_ZK_M(Y)\), and
\(P=S_u V\). Both \(V\) and \(P\) are unimodal by minimality
and convolution. Since the arms are noninterval, \(u\ge2\)
and \(b\ge3\). Write the three outer weights as \(a=R-\delta\), with \(\delta=1,2,3\).

Suppose first that \(Q_M(X,Y)\) is nondecreasing through \(R\).
At the next rank,
\(\Delta (Q_M)_{R+1}=V_R-V_{R-u}-1\) by
\eqref{eq:relative-append}. A failure then requires
\(V_R=V_{R-u}\) and a later strict increase of \(P\).
If a value strictly larger than the common endpoint
occurred inside \([R-u,R]\), unimodality of \(V\) would
give \(V_i-V_{i-u}\le0\) for all \(i>R\). Indeed the
entering values have passed the modal interval, whereas
an outgoing value before that interval is at least
\(V_{R-u}\) and any outgoing value after it is at least the
entering value. Thus a failure requires \(V\) to be
constant throughout \([R-u,R]\).

If \(2\le\delta\le u\) and \(b\ge\delta\), this constancy
is impossible because \(V_{R-\delta+1}>V_{R-u}\).
To prove the strict comparison, all interval products in
\(V=F_Z K_M(Y)\) have endpoint sum at least
\(\wt(Z)+R+1=2R-\delta-u+1\), the sum of the two indices.
By Lemma~\ref{lem:tail-budget}, the height interval of
\(K_M(Y)\) just above its coefficient at \(\delta\)
starts at \(l\ge\delta+1\) and contains \(R-\delta+1\).
Its convolution with the lowest interval of \(F_Z\)
has at index \(i=R-u,R-\delta+1\) the coefficient
\(\max(0,\min(\wt(Z),i-l)+1)\). The earlier index satisfies
\(i-l\le\wt(Z)-1\), the later index is at least \(l\), and
the indices are distinct. This layer contributes a strict increase.

For \(\delta=1\), growth through \(R\) follows from
Lemma~\ref{lem:kernel}. Each interval product in \(P\)
has center at least \(R+1/2\), so all contribute
nonnegatively to \(\Delta P_{R+1}\).
If \([q,R]\) is the lowest height interval of \(K_M(Y)\),
its product with \(S_u\) and the lowest interval of
\(F_Z\) contributes \(\coef{R-q}{[u][R-u]}\ge1\)
to this difference. Hence \(\Delta (Q_M)_{R+1}\ge0\), and from that rank onward the
coefficients agree with the unimodal window.
For \(\delta=2\) growth again reaches \(R\), and the
strict comparison above excludes constancy of \(V\) on \([R-u,R]\).

For \(\delta=3\), growth initially reaches \(R-1\).
Each interval product in \(P\) has center at least
\(R-1/2\). The height intervals of \(K_M(Y)\) of height
at least two that contain \(R-1\) have right endpoint
\(R-1\) and left endpoint at least three, because the
coefficients at \(2\) and \(R\) are at most one.
Their products with \(S_u\) and the lowest interval of
\(F_Z\) each contribute at least one to \(\Delta P_R\). Indeed, writing \(W=S_u[\,\wt(Z)+1\,]\) and \(l\) for the left
endpoint of the height interval, the difference is
\(W_{R-l}-W_0\). Here \(1\le R-l\le R-3=\wt(Z)+u\),
so \(W_{R-l}\ge1\) and \(W_0=0\). The \(K_M(Y)\) contribution is
\(1-(K_M(Y))_{R-1}\), which these intervals offset.
Thus \(Q_M\) increases through \(R\). For \(u\ge3\) the strict
comparison for \(\delta\le u\) rules out a plateau of \(V\).

For \(u=2\), put \(c=\wt(Z)=R-5\), \(A=F_Z\), and
\(B=F_Y\). If \(b\ge4\), we show \(\Delta V_{R-1}>0\).
At this edge \(\Delta V\) equals \(\Delta(BK_M(Z))\), since both
\(A\) and \(B\) vanish at the two ranks being compared.
For the lowest interval \([0,b]\) of \(B\), with
\(L=M-3\), the contribution is
\[
 \sum_{j=c-L+1}^c A_j-\sum_{j=0}^L A_j+A_{L+1}\ \ge0 .
\]
Indeed the higher-layer endpoint bound reflects
\(j=1,\ldots,L\) to \(c+1-j\) with at least as great a
coefficient, and \(A_{L+1}\ge A_0=1\).
For a higher interval \([l,r]\), its contribution is
\[
 (K_M(Z))_{b-l+M-1}-(K_M(Z))_{b-r+M-2}\ge0 .
\]
The indices are ordered and their sum is at most
\(b+2M-4=c+M+1\). For the second interval
\([2,b-1]\), the contribution is
\((K_M(Z))_{c+2}-(K_M(Z))_{M-1}\ge A_{M-1}\ge1\)
by \eqref{eq:reflected-lower} with \(v=M-2\).
Thus the required plateau is excluded.
If \(b=3\), put \(C=AB\), so \(\deg C=R-2\), \(C_1=2\),
and \(C_{R-2}=1\).
The relative identity gives
\(\Delta V_{R-1}=2C_{R-2}-C_{R-1}-C_1=0\) and
\(\Delta V_{R}=-C_2<0\). The strict fall at \(R\) rules out the required plateau in this last case.
\end{proof}
\section{Estimates for the mixed case and its finite verification}\label{app:finite}

For the mixed comparison, write \(A=F_X\), \(B=F_Y\), \(a=\wt(X)\),
and \(b=\wt(Y)\). At the reverse index \(n\) put \(j=a+b-n\) and
\(r=a-j\). Suppose \(B_{j+1}>B_j\) and \(A_{j+1}\le A_j\), and let
\(K=\kappa(Y)\) be the first weight missing from the ideal of \(Y\).
The first few nonzero increments of \(\mathcal E_Y\) already give
a lower bound for the sum subtracted in \eqref{eq:mixed}. For
\(K\le12\) this bound leaves only finitely many possibilities for
the other arm \(X\) and the rank \(r\).

\begin{lemma}[Early source increments]\label{lem:early-source}
If a nonempty increasing arm \(Y\) has first missing weight \(K\),
then the forward source satisfies
\[
 (\mathcal E_Y)_K\ge1,\qquad
 (\Delta\mathcal E_Y)_{K+i}\ge2^{i-2}
 \quad(2\le i\le K-3).
\]
\end{lemma}
\begin{proof}
Normalize the increasing arm and write \(J=z^{-K}\mathcal E_Y\).
We first bound its coefficients by \(J_0\ge1\) and \(J_i\ge2^{i-1}\)
for \(1\le i\le K-2\). A string of ones has \(K=2\). If appending \(t\)
raises the first missing weight from \(k\) to \(k+1\), the normalized
source changes from \(J\) to \([t]J+z^{t-k}\), with \(t\ge k\).
At the required rank \(i\) the window contains \(J_0,\ldots,J_{i-1}\),
whose lower bounds sum to \(2^{i-1}\). If the first missing weight
stays at \(k\), then \(t=k-1\) and the new source is \(1+z[t]J\).
The same prefix gives the same lower bound.

Now take first differences at \(2\le i\le K-3\). In an append preserving
the first missing weight, the new difference at \(i\) is the old
coefficient at \(i-1\), at least \(2^{i-2}\). In an append raising the
first missing weight, the departing coefficient vanishes because \(i<t\).
The extra monomial can lower the difference by at most one, leaving at
least \(2^{i-1}-1\ge2^{i-2}\).
\end{proof}

\subsection{An interval product and a finite degree bound}

All source increments through \(j+1\) in \eqref{eq:mixed} are
nonnegative. In reverse order \(A_{j+1-K-i}=q_{X,r+K-1+i}\),
since \(r=a-j\) and both sequences are extended by zero. Lemma~\ref{lem:early-source}
therefore gives
\begin{equation}\label{eq:finite-payment}
 T\le C_0-q_{X,r+K-1}-\sum_{i=2}^{K-3}2^{i-2}q_{X,r+K-1+i}.
\end{equation}
Here \(C_0=\Delta q_{X,r+1}-\Delta q_{X,r}\). When \(r\ge0\) and
\(T>0\), both slopes are positive.

Write the normalized first arm as \(X=1^\alpha(t_1,\ldots,t_h)\)
and put \(D=\sum_i(t_i-1)\) and \(K_X=\prod_i[t_i]\).
Then \(a=\alpha+h+D\). To arrange that \(t_1-1\le\alpha\), it suffices,
when \(t_1-1>\alpha\), to
exchange \([\alpha+1]\) and \([t_1]\) in
\(F_{1^\alpha t_1}=1+z[\alpha+1][t_1]\). The new initial run has
length \(t_1-1\), and the next part is \(\alpha+1\le t_1\). The
arm remains increasing and preserves its weight and rank polynomial.

\begin{lemma}[Comparison with the interval product]\label{lem:skeleton}
The sequence \(v_i=(\Delta q_{X,i})_+\) vanishes outside
\(0\le i\le D\) and satisfies \(v_i(K_X)_j\ge v_j(K_X)_i\) for \(i<j\).
For \(0\le r\) with \(v_{r+1}>v_r>0\), let \(e=D-2r-2\).
Then \(e\ge0\) and
\begin{equation}\label{eq:skeleton-curvature}
 (r+1)(v_{r+1}-v_r)\le(e+1)v_r .
\end{equation}
\end{lemma}
\begin{proof}
For an all-ones arm, \(v\) and \(K_X\) are both the unit sequence.
In an append \(X=Zt\), write \(W^+=[t](\Delta q_Z)_+\) and
\(W^-=[t](-\Delta q_Z)_+\). The positive differences of \(q_Z\)
precede its negative differences by unimodality. Thus \(W^+,W^-\)
are in ratio order by Section~\ref{sec:transport}, and
\((W^+-W^-)_+/W^+\) is nonincreasing where defined. The same
convolution rule applied to the induction hypothesis makes
\(W^+/K_X\) nonincreasing where defined. Multiplying gives the
desired comparison for \(v'=(W^+-W^-)_+\). In the cross-product
inequality with \(i<j\), if \(v'_j(K_X)_i>0\), then \(W^+_j>0\),
and the ratio order of \(W^+,K_X\) implies \(W^+_i>0\) and
\((K_X)_j>0\). These are all the denominators needed for the two ratios.
If \(v'_j(K_X)_i=0\), nonnegativity suffices.
The terminal unit \(z^{\wt(Z)+t}\) changes a final difference
of \(-1\) to zero and adds a subsequent negative difference.
Thus \(v=v'\), and the comparison at \(i=0\) gives \(v_j=0\)
for \(j>D\), since \((K_X)_0=1\) and \(K_X\) has degree \(D\).

For an interval product \(J=\prod_i[d_i+1]\) of degree
\(D=\sum_i d_i\), and \(2r\le D-1\), one has
\begin{equation}\label{eq:degree-ratio}
 (r+1)J_{r+1}\le (D-r)J_r.
\end{equation}
The right side minus the left is the coefficient of \(z^r\) in
\(DJ-(1+z)J'\). In its product-rule expansion, the factor
\(d_i[d_i+1]-(1+z)([d_i+1])'\) has coefficients
\(d_i-2u-1\) for \(0\le u<d_i\), with the last coefficient zero.
Pair \(u\) with \(d_i-1-u\), taking \(u\) on the lower side.
The remaining interval product is symmetric and unimodal
of degree \(D-d_i\). Its coefficient at \(r-u\) is at
least that at \(r-(d_i-1-u)\), since their sum is at
most \(D-d_i\). Thus every pair contributes nonnegatively to \eqref{eq:degree-ratio}.
If \(v_{r+1}>v_r>0\), the ratio order forces
\((K_X)_{r+1}>(K_X)_r\), giving \(D\ge2r+2\) by symmetry.
Then \eqref{eq:degree-ratio} and the ratio order give
\((r+1)v_{r+1}\le(D-r)v_r\), yielding
\eqref{eq:skeleton-curvature}.
\end{proof}

\begin{lemma}[A finite degree reduction]\label{lem:owner-bound}
In the mixed range with \(r\ge0\), \(h\ge2\) and \(T>0\),
put \(d=K-\alpha\). Then
\begin{equation}\label{eq:finite-domain}
 0\le e=D-2r-2\le d-3,\qquad
 D\le \min(2d(d-3),\,12d-25).
\end{equation}
\end{lemma}
\begin{proof}
Since \(A_{j+1}\le A_j\),
\eqref{eq:mixed} gives \(v_{r+1}>v_r>0\), using
initial support. Lemma~\ref{lem:skeleton} gives \(e\ge0\).
If \(\deg J_X=\alpha+D\ge2r+K\), symmetry and unimodality give
\((J_X)_{r+K-1}\ge(J_X)_{r+1}\). Lemma~\ref{lem:head-range}
then gives \(q_{X,r+K-1}\ge(J_X)_{r+1}\ge C_0\).
The first term subtracted in \eqref{eq:finite-payment} would then
imply \(T\le0\). Thus \(e\le d-3\). Put \(L=d-e-3\).
Each shifted interval product in \(q_X\) has endpoint sum at least
\(\alpha+D\), as in Lemma~\ref{lem:head-range}. For \(r\ge L\),
reflection therefore gives \(q_{X,r+K-1}\ge q_{X,r-L}\), since
\[
 (r-L)+(r+K-1)=\alpha+D .
\]

For the quadratic bound in \eqref{eq:finite-domain}, suppose \(r\ge dL+e\).
By the degree ratio and Lemma~\ref{lem:skeleton},
\[
 v_r\le P v_{r-L},\quad
 P=\prod_{i=0}^{L-1}(1+x_i),\qquad
 x_i=\frac{e+3+2i}{r-i}.
\]
If \(L=0\) the product is one.
With \(B=r-L+1>0\), the sum \(S=\sum_i x_i\) satisfies
\(BS\le L(e+L+2)=L(d-1)\).
Consequently \(B(1-S)\ge r+1-dL\ge e+1>0\).
For nonnegative \(x_i\) of sum \(S\le1\), we have
\((1-S)\prod_i(1+x_i)\le1\), by iteration of
\((1-u-x)(1+x)\le1-u\) for a partial sum \(u\).
Hence \((e+1)P\le r+1\), and \eqref{eq:skeleton-curvature} gives
\(C_0\le v_{r-L}\le q_{X,r-L}\), so the first source
term is at least \(C_0\). Thus \(r<dL+e\).
As all parameters are integral,
\(D=2r+2+e\le2d(d-3)-(2d-3)e\le2d(d-3)\).

For the linear bound suppose \(D\ge12d-24\).
Substitute \(d=e+L+3\) to get
\(r\ge3e+6L+5\) and \(r\ge4e+3\).
When \(L=0\), the first source coefficient suffices since \(r\ge e\).
When \(L=1\) the two bounds are
\(r v_r\le (r+e+3)v_{r-1}\) and
\((e+1)(r+e+3)\le r(r+1)\). The latter follows at
\(r=2e+3\) by expansion, and its difference increases
thereafter. In both cases the first source coefficient
dominates \(C_0\), as above.

For \(L\ge2\), \(L\le K-3\), and the source term at index
\(K+L\) is at least \(2^{L-2}q_{X,r+K-1+L}\) by
Lemma~\ref{lem:early-source}. The forward index \(j+1-K-L=r+h-2L\) is at least \(h\),
so this term occurs in the sum.
The reflected index is \(r-2L\), making this source term
at least \(2^{L-2}v_{r-2L}\). The backward ratios satisfy
\[
 v_r\le v_{r-2L}\prod_{i=0}^{2L-1}\frac{r+e+3+i}{r-i}.
\]
Pair the factors with indices \(i\) and \(2L-1-i\). In each pair the
numerator and denominator are obtained by adding \(i(2L-1-i)\) to
\((r+e+3)(r+e+2L+2)\) and \(r(r-2L+1)\), respectively. The first of
these products is at most twice the second. Indeed, at
\(r=3e+6L+5+x\), twice the second minus the first is
\(2e^2+4eL+6e+6L+4+(4e+6L+7)x+x^2\ge0\).
Each paired ratio is therefore at most two, and their product is at
most \(2^L\).

Since \(r+1\ge4(e+1)\),
\eqref{eq:skeleton-curvature} gives
\(C_0\le2^{L-2}v_{r-2L}\). This again contradicts
\(T>0\), proving the linear bound.
\end{proof}

\subsection{Exact enumeration of the residual domain}

\begin{proposition}[Small first missing weights]\label{prop:small-k}
Under the hypotheses of Condition~\ref{cond:local}, suppose
\(A_{j+1}\le A_j\) and \(B_{j+1}>B_j\) at \(j=a+b-n\).
If \(K=\kappa(Y)\le12\), then \eqref{eq:local} holds.
\end{proposition}
\begin{proof}
By Proposition~\ref{prop:mixed-head} we may take \(0\le r<k\),
\(j-2>s_Y\), and \(h\ge2\). Lemma~\ref{lem:capacity} settles every
arm \(X\) whose largest part is at least \(K-1\), hence all of \(K=3\).
Suppose \(T>0\) at a remaining rank. In Lemma~\ref{lem:owner-bound},
\(D\ge h\ge2\) forces \(d\ge4\), since \(e\ge0\) gives \(d\ge3\)
and the quadratic bound gives \(D\le0\) at \(d=3\).
The increments \(t_i-1\) form a partition of \(D\) with parts in
\([1,K-3]\) and smallest part at most \(\alpha\). Also
\(j-2>s_Y\ge S_K\), so the number \(h\) of parts satisfies
\begin{equation}\label{eq:height-small}
 \max(2,S_K+3-\alpha-D+r)\le h\le D .
\end{equation}

For a rank state \(K,\alpha,r,e\) with
\(D=2r+2+e\) and \eqref{eq:finite-domain}, set
\(L=d-e-3\) and
\[
 b_m=\begin{cases}
 \binom r m/\binom{r+e+m+2}{m},&0\le m\le r,\\
 0,&m>r .
 \end{cases}
\]
The degree ratio \eqref{eq:degree-ratio}, applied
successively with the ratio order of Lemma~\ref{lem:skeleton}, gives
\(v_{r-m}\ge b_m v_r\) for \(m\le r\).
For every used source offset \(i\), reflection at
endpoint sum \(\alpha+D\) gives
\(q_{X,r+K-1+i}\ge q_{X,r-L-i}\) when \(L+i\le r\).
An offset with \(L+i>r\) simply has the zero lower bound
\(b_{L+i}v_r\).
Consequently the source in \eqref{eq:finite-payment}
is at least \(C_0\) by \eqref{eq:skeleton-curvature} whenever
\begin{equation}\label{eq:finite-ratio}
 (r+1)\left(b_L+\sum_{i=2}^{K-3}2^{i-2}b_{L+i}\right)\ge e+1.
\end{equation}
The calculation starts with all integer tuples \((K,\alpha,r,e)\)
satisfying \eqref{eq:finite-domain}, with \(4\le K\le12\),
\(1\le\alpha\le K-4\), and \(r\ge0\). Of the \(4\,169\) tuples, \(641\) admit no
length \(h\) in \eqref{eq:height-small} and another \(2\,852\)
satisfy \eqref{eq:finite-ratio}. For the remaining \(676\), we
enumerate partitions with the indicated part and length bounds and
retain the arm and rank pairs with \(v_{r+1}>v_r>0\).
Table~\ref{tab:finite} gives the counts by \(K\).

\begin{table}[htbp]
\centering
\caption{Finite verification by first missing weight \(K\)}
\label{tab:finite}
\small
\begin{tabular*}{0.93\linewidth}{@{\extracolsep{\fill}}rrrrrrr@{}}
 \toprule
 \(K\) & States & Height & Ratio & Remaining & Arm--rank inputs & Critical\\
 \midrule
 4 & 0 & 0 & 0 & 0 & 0 & 0\\
 5 & 7 & 1 & 4 & 2 & 4 & 3\\
 6 & 35 & 3 & 22 & 10 & 45 & 31\\
 7 & 103 & 12 & 69 & 22 & 289 & 145\\
 8 & 237 & 28 & 166 & 43 & 1\,488 & 569\\
 9 & 441 & 53 & 311 & 77 & 7\,667 & 2\,373\\
 10 & 719 & 99 & 500 & 120 & 38\,053 & 7\,349\\
 11 & 1\,083 & 172 & 742 & 169 & 148\,639 & 19\,591\\
 12 & 1\,544 & 273 & 1\,038 & 233 & 602\,239 & 49\,769\\
 \midrule
 Total & 4\,169 & 641 & 2\,852 & 676 & 798\,424 & 79\,830
\\
 \bottomrule
\end{tabular*}\par\smallskip
\begin{minipage}{0.93\linewidth}
\footnotesize\singlespacing\raggedright
\textit{Notes.} Height counts states with no admissible length in
\eqref{eq:height-small}. Ratio counts those settled by
\eqref{eq:finite-ratio}. The remaining states determine the arm--rank
inputs. Critical inputs have positive reverse slope and positive curvature.
\end{minipage}
\end{table}

For each retained pair we check
\[
 \Delta q_{X,r+1}-\Delta q_{X,r}
 \le q_{X,r+K-1}+\sum_{i=2}^{K-3}2^{i-2}q_{X,r+K-1+i}.
\]
One program constructs the reversed coefficients by appending the
parts of \(X\) and indexes the rank states by \(e,r\). A second program
indexes them by \(D,r\), generates partitions by their largest part,
and obtains the forward coefficients by summing over word lengths
as in Section~\ref{sec:transport}. The programs give the same
\(79\,830\) critical rows, with agreement on the arm, the rank,
and both sides of the inequality. The smallest value of the right side
minus the left is \(2\). Thus \eqref{eq:finite-payment} gives
\(T\le0\) throughout the remaining range. For each such \(X,r\) the
bound applies to every admissible \(Y\) by Lemma~\ref{lem:early-source}.
\end{proof}

\subsection{Proof of the slope bounds for large first missing weight}
\label{app:slope}

\begin{proof}[Proof of Lemma~\ref{lem:dyadic-slope}]
Put
\(\rho=3/5\), \(\gamma=35/32\), \(\eta=125/22\),
\(\theta=\gamma^2/(2\rho)=6125/6144\),
\(v=\gamma/(2\rho)=175/192\).
Consider any arm \(Z\) formed by appending parts of size at least two
to an initial run of \(\beta\ge0\) ones. Let \(F=F_Z\) and let \(M\)
be its main interval degree, the sum of \(\beta\) and \(t-1\) over the
appended parts \(t\). The append recurrence gives
\begin{equation}\label{eq:eval-bound}
 F(\rho)+\eta\le10\gamma^M.
\end{equation}
Indeed the initial interval \([\beta+1]\) evaluates at
most \(5/2\) for \(\beta\ge0\).
The factors \(s_t=S_t(\rho)\) obey
\(s_t\le\gamma^{t-1}\) for \(t\ge2\). Direct substitution covers \(2\le t\le5\), and
\(s_t\le3/2\le\gamma^5\) thereafter.
Write \(F_{\rm app}=1+S_tF\) for the polynomial after an append.
For \(t\ge3\), \(\eta(s_t-1)\ge1\), and hence
\(F_{\rm app}(\rho)+\eta\le s_t(F(\rho)+\eta)\).
For \(t=2\) the remaining margin
\(10(\gamma-24/25)-1-\eta/25=97/880>0\)
proves the induction. Put \(r_t=s_t/\gamma^{t-1}\).
The ratios \(r_t\) decrease for \(t\ge4\) since
\(\rho^{t+1}\le(\gamma-1)s_t\) there, and \(r_3\le r_4\).

First suppose \(2m\ge M_Y\).
If \(K\ge13\), the last eleven parts dominate \(2,\ldots,12\)
componentwise. Remove the last eight parts. Let \(G\) be the rank
polynomial of the remaining prefix and \(M_G\) the degree of its main
interval product. The last two parts in this prefix are at least three
and four, respectively, so
\(G(\rho)\le10\gamma^{M_G}r_4^2\).
The eight removed factors have total degree at least \(60\).
Let \(Y^-\) be \(Y\) with its last part \(t\) removed. By
\eqref{eq:append}, for \(m\ge2\),
\(\Delta B_m=\coef{m-1}{F_{Y^-}}-\coef{m-1-t}{F_{Y^-}}\).
In particular \(\Delta B_m\le\coef{m-1}{F_{Y^-}}\).
Expanding its last seven append steps,
\[
 F_{Y^-}=\sum_{\ell=0}^6\prod_{\text{last }\ell}S_t
              +\left(\prod_{\text{last }7}S_t\right)G .
\]
The seven proper suffix terms count words with at most six parts. At
weight \(m-1\), their at most five cuts lie in \(1,\ldots,m-2\).
Recording the cuts in order and padding with zeros gives at most
\(m^5\) such words.
For the product of seven unshifted intervals,
\(\coef i{\prod[t]}\rho^i\le\binom{i+6}{6}\rho^i
\le3003\rho^8\).
The last bound follows from the adjacent ratio
\(\rho(i+7)/(i+1)\), with maxima at eight and nine.
After the seven shifts and extraction at rank \(m-1\), the last
term of \(F_{Y^-}\) contributes at most \(3003\rho^{16}G(\rho)\)
when multiplied by \(\rho^m\).
Since \(M_Y\ge67\), we have \(m\ge34\). The constants satisfy
\[
 25\cdot30030\rho^{16}r_4^2\le\gamma^{60},
 \qquad 100\cdot34^5\le2^{34},\qquad (35/34)^5<2 .
\]
Together with \(M_G+60\le2m\) and \(\theta^m\le1\),
these inequalities bound the contribution of \(G\), divided by \(2^m\),
by \(1/25\). The corresponding bound for the proper suffixes is
\(1/100\), since \(m^5/2^m\) decreases for \(m\ge34\). Their sum is
\(1/20\), as required in \eqref{eq:constant-density}.

For \eqref{eq:rank-density}, the last fourteen parts dominate
\(2,\ldots,15\) componentwise.
Let \(w,u\) denote the two parts immediately before the last one,
and let \(G\) be the rank polynomial of the prefix before \(w\). Its evaluation is bounded by
\(10\gamma^{M_G}P_*\), where
\(P_* = r_4^2\prod_{t=5}^{12}r_t\). This follows by including
the first of the fourteen in the prefix and applying the offset inequality
to the next ten.
The removed main degree is at least \(39\).
The prefix before the last two parts has main degree at
least \(79\). The last part is at least \(u\), so \(M_Y\ge2u+4\) and \(m-1>u\).
Hence in \(F_{Y^-}=1+S_u+S_uS_wG\) the first two
terms do not contribute at \(m-1\). The coefficient bound
\((i+1)\rho^i\le6/5\) for \([u][w]\), obtained by
consecutive ratios, gives
\[
 \gamma^{39}\frac{\Delta B_m}{2^m}
 \le10\rho^3(6/5)P_*\theta^m.
\]
Here \(\theta^{16}\le20/21\) and
\(20(20/21)^{19}\le8\).
Writing \(m=16a+r\), \(0\le r<16\), the adjacent ratios
of \((a+1)(20/21)^a\) give \(m\theta^m\le128\).
The exact inequality
\((41472/125)P_*\le\gamma^{39}\) now proves
\eqref{eq:rank-density}.

Now put \(O=M_Y-B_Y\) and suppose \(m>O\). If the last part gives
a largest interval factor, deleting that factor leaves
degree \(O\) and gives \(\Delta B_m/2^m\le5v^{m-1}\)
by \eqref{eq:eval-bound}. If the initial interval
\([\beta+1]\) is largest, expand
\(F_Y=A+\left(\prod S_t\right)[\beta+1]\) with \(A\)
the proper suffix sum. Then
\((1-z)F_Y\) is bounded coefficientwise by
\(A+\prod S_t\), the formal zero-initial-run polynomial
of degree \(O\) in the main interval product. This gives
\(\Delta B_m/2^m\le10v^m\).
For \(K\ge13\), the respective lower bounds on \(O\)
are \(56\) and \(66\). Use \(100v^{56}\le1\) and
\(200v^{67}\le1\).
For \(K\ge16\), \(O\ge92\) in either case.
Since \(10\cdot81v^{80}\le1\) and \(10(i+1)v^i\) decreases for
\(i\ge80\), we obtain \eqref{eq:rank-density}.
Since \(s_Y=\min(\lfloor M_Y/2\rfloor,M_Y-B_Y)\), the two ranges
cover every integer \(m>s_Y\).
\end{proof}
\end{document}